\documentclass[a4paper,10pt]{article}

\usepackage{amsmath}
\usepackage{amssymb}
\usepackage{amsthm}
\usepackage{amsfonts}
\usepackage{xcolor}
\usepackage{graphicx}
\usepackage{mathtools}
\usepackage{enumerate}
\usepackage{multirow}
\usepackage{rotating}
\usepackage{setspace}
\usepackage{lineno}
\usepackage{url}

\theoremstyle{theorem}
\newtheorem{thm}{Theorem}

\newtheorem{lem}{Lemma}

\theoremstyle{definition}

\usepackage{geometry}

\begin{document}
	
	\title{ \large{A Nested Multi-Scale Model for COVID-19 Viral Infection }}
	
	\vspace{0.1in}
	%{\bf{\Large }}
	\author{{\small Bishal Chhetri$^{a,*}$, D. K. K. Vamsi$^{a}$,  Carani B Sanjeevi$^{b}$}\hspace{2mm} \\
		{\it\small $^{a}$Department of Mathematics and Computer Science, Sri Sathya Sai Institute of Higher Learning, Prasanthi Nilayam}, \\
		{\it \small Puttaparthi, Anantapur District - 515134, Andhra Pradesh, India}\\
		%{\it\small $^{b}$Central Leprosy Teaching and Research Institute - CLTRI, Chennai, India}\\
		%{\it\small $^{c}$ Professor and Head, Department of Community Medicine, All India Institute of Medical Sciences - AIIMS, Nagpur,India}\\
		{\it\small $^{b}$  Vice-Chancellor, Sri Sathya Sai Institute of Higher Learning -  SSSIHL, India}\\
		%{\it\small $^{c}$ Department of Medicine, Karolinska Institute, Stockholm, Sweden }\\
		{\it\small bishalchhetri@sssihl.edu.in, dkkvamsi@sssihl.edu.in$^{*}$,}\\
		%{\it\small ananthvs@sssihl.edu.in, prakashdmacs@gmail.com,      swapnamuthuswamy@gmail.com,}\\
		{\it\small sanjeevi.carani@sssihl.edu.in, sanjeevi.carani@ki.se}\\
		%{\small $^{1}$  First Author},
		{ \small $^{*}$ Corresponding Author}
		\vspace{1mm}
	}
	
	\date{}
	\maketitle

\begin{abstract} \vspace{.25cm}

 In this study, we develop and analyze a nested multi-scale model for COVID -19 disease that integrates  within-host scale and between-host scale sub-models. First, the well-posedness of the multi-scale model is discussed, followed by the stability analysis of the equilibrium points. The disease-free equilibrium point is shown to be globally asymptotically stable for $R_0 < 1$. When $R_0$ exceeds unity, a unique infected equilibrium exists, and the system is found to undergo a forward (trans-critical) bifurcation at $R_0=1$. Two parameter heat plots are also done to find the parameter combinations for which the equilibrium points are stable. The parameters $\beta, \pi$ and $\Lambda$ are found to be most sensitive to $R_0$. The influence of within-host sub-model parameter on the between-host sub-model variables is numerically illustrated. The spread of infection in a community is shown to be influenced by within-host level sub-model parameters, such as the production of viral particles by infected cells $(\alpha)$, the clearance rate of infected cells by the immune system $(x)$, and the clearance rate of viral particles by the immune system $(y)$. The comparative effectiveness of the three health interventions (antiviral drugs, immunomodulators, and generalized social distancing) for COVID-19 infection  was examined using the effective reproductive number $R_E$ as an indicator of the effectiveness of the interventions. The results suggest that a combined strategy of antiviral drugs, immunomodulators and generalized social distancing would be the best strategy to implement to contain the spread of infection in the community. We believe that the results presented in this study will help physicians, medical professionals, and researchers to make informed decisions about COVID -19 disease prevention and treatment interventions.
 
\end{abstract}

\section*{keywords}
Multi-Scale Models, Well Posedness, Elasticity, Comparative Effectiveness, Heat Plot

\section{Introduction}

COVID -19 is a contagious respiratory and vascular disease that has resulted in more than 209 million cases and 4.4 million deaths worldwide. It is caused by severe acute respiratory syndrome coronavirus 2 (SARS-CoV-2). It was declared an international public health emergency on January 30. 

Mathematical modeling of infectious diseases is one of the  highly researched area in applied mathematics. Mathematical epidemiology has contributed to a better understanding of the dynamic behavior of infectious diseases, their impact, and possible future predictions of their spread. In general, the transmission of an infectious disease system is a consequence of processes acting on at least two different
scales, ie, within-host and between-host scales \cite{handel2015crossing, gog2015seven}. At within-host scale there are three ways in which the transmission of an infectious disease system is regulated that is, when the immune system impacts the dynamics of an infectious disease system in a population: first, by regulating pathogen dynamics at within-host scale; second, by affecting herd immunity at between-host scale; and finally, by altered immuno-demography through change in the immune status due to age, disease, or therapy. There are number of studies on infectious disease systems which has established that transmission of infectious diseases depends critically on within-host scale processes. In particular, these studies established that the transmission potential of an infectious host increases with increasing pathogen load in the infected host \cite{handel2015crossing}. Studies on HIV and dengue virus transmission established a sigmoid functional relationship between host infectiousness and within-host scale pathogen load \cite{mellors1996prognosis, quinn2000viral}. Other studies that have also confirmed a similar functional relationship between host infectiousness and within-host scale pathogen load include\cite{jeffery1955infectivity} for malaria, and \cite{kaplan1996male} for human
T lymphotropic virus.  Studies such as \cite{feng2012model, almocera2018multiscale} considered transmission rate at between-host scale as a linear function of viral load. Several within-host and between-host compartment models are developed to study the spread of COVID-19 infection at individual level as well as at community level. Some of the important mathematical modelling studies that deal with transmission and spread of COVID-19 at between-host level can be found in  \cite{chen2020mathematical, 7, 4, 6,samui2020mathematical, ndairou2020mathematical, zeb2020mathematical,leontitsis2021seahir, wang2020four, dashtbali2021compartmental, zhao2020five, biswas2020covid, sarkar2020modeling}. Within-host mathematical modelling studies that deal with the interplay between the viral dynamics and immune response of COVID-19 disease can be found in \cite{chhetri2021within, vargas2020host}.

A thorough understanding of infectious disease transmission requires knowledge of the processes at different levels of infectious disease and of the interplay between these levels. To get a clear idea of the dynamics of the disease, the  different time scale models need to be integrated. Multi-scale models of infectious disease systems integrate the within-host scale and the between-host scale. There are five main different categories of multi-scale models that can be developed at the different levels of organization of an infectious disease system, which are: individual-based multiscale
models (IMSMs), nested multiscale models (NMSMs), embedded multiscale models (EMSMs), hybrid multiscale
models (HMSMs), and coupled multiscale models (CMSMs) \cite{garira2017complete}. The multi-scale models for influenza infections can be found in \cite{heldt2013multiscale, guo2015multi}. 

Multi-scale modeling of COVID -19 disease is still in its infancy, and to our knowledge there is only one literature in which a multi-scale model has been developed for COVID -19 disease \cite{prakash2020control}. A multi-scale model would be extremely helpful in understanding the spread of COVID -19 infection and evaluate the efficacy of the interventions not only at individual level but also the at the population level. Therefore, in this article, we develop a nested multi-scale model that integrates within-host and between-host sub-models. As the importance of multi-scale modeling in disease dynamics is increasingly recognized, we believe that our study contributes to the growing  knowledge on multi-scale modeling of COVID -19 disease. The work in this study is aimed at two types of audiences - disease modelers and public health planners. For disease modelers, this study provides an alternative approach to modeling acute viral infections. For public health planners, this study provides a model-based approach to evaluating the effectiveness of public health interventions.  

The paper is organised as follows: In section 2, we develop a nested multi-scale model, study the stability and bifurcation analysis and numerically illustrate the theoretical results obtained. We also do the sensitivity of $R_0$ with the model parameters and heat plots. In section 3 we evaluate the comparative effectiveness of the three health interventions. The study concludes with section 4 where we summarize the findings of the work.

\section{Multi-Scale Model}
The Multi-Scale model that we develop describes  COVID-19 disease dynamics across two different time scales, that is, within-host scale and  between-host scale. The 
 multi-scale model is based on the monitoring of seven variables, namely, susceptible epithelial cells $U$, infected epithelial cells $U^*$, viral load $V$, and four between-host variables namely,  susceptible human $S$, infection human population $I$, exposed human population $E$, and recovered human population $R$. We make the following assumptions for the multi-scale 
 model. \\
 
 a) The dynamics of the within-host scale variables are assumed to occur at fast time scale $s$ so that $U = U(s)$, $U^* = U^*(s)$ and $V = V(s)$ while the dynamics of the between-host scale variables is assumed to occur at slow time scale $t$ so that $S = S(t)$, $I = I(t)$, $E=E(t)$, and $R=R(t)$.\\
 
 b) The transmission rate $\beta$ and disease induced death rate $d$ at between-host scale are assumed to be a functions of viral load.\\
 
 c) Immune response is captured in the model through the parameters $b_1, b_2, b_3, b_4, b_5, b_6$ and $d_1, d_2, d_3, d_4, d_5, d_6$  where these parameters denote the rate at which infected cell and virus are cleared by the release of cytokines and chemokines such as IL-6 TNF-$\alpha$, INF-$\alpha$, CCL5, CXCL8 , and CXCL10.
 
 Based on the assumptions above the multi-scale model is described by the following system of differential equations.

\begin{eqnarray}
   	\frac{dU}{ds}& =&  \omega \ - k U(s)V(s)  - \mu_{c} U(s)  \label{sec2equ1} \\
   	\frac{dU^*}{ds} &=& k U(s)V(s) \ -  { \bigg(d_{1}  + d_{2}  + d_{3} +  d_{4} + d_{5}+ d_{6}\bigg)U^*(s) }  \ - \mu_{c} U^*(s)   \label{sec2equ2}\\ 
   	\frac{dV}{ds} &=&  \alpha U^*(s)   \ -  \bigg( b_{1}  + b_{2}  + b_{3}  +  b_{4}+ b_{5} + b_{6}\bigg)V(s)    \ -  \mu_{v} V(s) \label{sec2equ3}\\
   	\frac{dS}{dt}& =&  \Lambda \ - \beta(V(s)) S(t)I(t)  - \mu S(t)  \label{sec2equ1} \\
   	\frac{dE}{dt} &=& \beta(V(s)) S(t)I(t) \ - (\mu+\pi+\gamma_1)E(t)    \label{sec2equ2}\\ 
   	\frac{dI}{dt} &=& \pi E(t)-(\mu+\gamma_2)I(t)-d(V(s))I(t)\label{sec2equ3} \\
   	\frac{dR}{dt}&=& \gamma_1 E(t)+ \gamma_2 I(t)-\mu R(t)
   \end{eqnarray} 
   
The meaning of each of the variables and parameters of the model is given in table 1.   \\

   \begin{table}[ht!]
     	\caption{Meanings of the Parameters}
     	\centering % used for centering table
     	\begin{tabular}{|l|l|} % centered columns (4 columns)
     		\hline\hline
     		
     		\textbf{Parameters} &  \textbf{Biological Meaning} \\  % inserts table
     		%heading
     		\hline\hline % inserts single horizontal line
     	  $\omega$ & natural birth rate of cells \\
     	\hline\hline
     	$\Lambda$ & birth rate of human population \\
     	\hline\hline
     		$\alpha$ & burst rate of the virus \\
     		
     	%	\hline\hline
     	%	$\alpha$ & rate of increase in viral load due to adverse effect \\
     		\hline\hline
     		$\mu$ & natural death rate of human population \\
     		\hline\hline
     		$\mu_{c}$ & natural death rate of cells \\
     		\hline\hline
     	
     		$\mu_v$ & natural death rate of virus \\
     		\hline\hline
     			$\pi$ & infection rate of exposed population  \\
     			\hline\hline
     				$k$ & infection rate of susceptible cell  \\
     				\hline\hline
     					$\gamma_1, \gamma_2$ & recovary rate of the exposed and infected human population  \\
     				\hline\hline
     	\end{tabular}
     \end{table} \vspace{.25cm}
 
 Considering experimental observations of the impact of viral load on disease transmission \cite{kawasuji2020transmissibility} and disease induced deaths, the linking of within-host and between-host sub-models is done considering $\beta=\beta(V(s))$ and $d=d(V(s))$. Though the exact functional relationship between viral load,  transmission rate and disease induced death rate is not known \cite{almocera2018multiscale}, as in \cite{feng2012model, almocera2018multiscale} we considered a linear form of the coupling functions $\beta$ and $d$: $\beta(V(s))=\beta V(s)$ and $d(V(s))=d V(s)$. We also observe that the first three equation of the between-host SEIR sub-model is independent of recovered population $R(t)$. Therefore without loss of generality we omit the last equation of the between-host SEIR sub-model. Taking $\beta$ and $d$ as a linear function of viral load and omitting the equation for $R(t)$ in between-host SEIR sub-model the final multi-scale model for COVID-19 disease dynamics is given by the following set of equations:

\begin{eqnarray}
   	\frac{dU}{ds}& =&  \omega \ - k U(s)V(s)  - \mu_{c} U(s)  \label{sec2equ1} \\
   	\frac{dU^*}{ds} &=& k U(s)V(s) \ -  { \bigg(d_{1}  + d_{2}  + d_{3} +  d_{4} + d_{5}+ d_{6}\bigg)U^*(s) }  \ - \mu_{c} U^*(s)   \label{sec2equ2}\\ 
   	\frac{dV}{ds} &=&  \alpha U^*(s)   \ -  \bigg( b_{1}  + b_{2}  + b_{3}  +  b_{4}+ b_{5} + b_{6}\bigg)V(s)    \ -  \mu_{v} V(s) \label{sec2equ3}\\
   	\frac{dS}{dt}& =&  \Lambda \ - \beta V(s) S(t)I(t)  - \mu S(t)  \label{sec2equ1} \\
   	\frac{dE}{dt} &=& \beta V(s) S(t)I(t) \ - (\mu+\pi+\gamma_1)E(t)    \label{sec2equ2}\\ 
   	\frac{dI}{dt} &=& \pi E(t)-(\mu+\gamma_2)I(t)-d V(s)I(t)\label{sec2equ3} 
   \end{eqnarray} 
 
 \subsection{\textbf{Reduced Multi-Scale Model for COVID-19 Dynamics}}
 
 The multi-scale model given by equations $(2.8)-(2.13)$ is not easy to analyze. The difficulty arises from the mismatch of time scales, since the within-host sub-model is in the form of a fast time scale $s$, while the between-host sub-model is in the form of a slow time scale $t$.
  To overcome these problems, we simplify the multi-scale model by using the within-host scale sub-model to define another quantity, which we use as a proxy for the infectivity of the individual host and which is called the area under the viral load curve \cite{hadjichrysanthou2016understanding, garira2020development}. Consider the within-host scale sub-model from the multi-scale model $(2.8)-(2.13)$.
  
  \begin{eqnarray}
   	\frac{dU}{ds}& =&  \omega \ - k U(s)V(s)  - \mu_{c} U(s)  \label{sec2equ1} \\
   	\frac{dU^*}{ds} &=& k U(s)V(s) \ -  { x U^*(s) }  \ - \mu_{c} U^*(s)   \label{sec2equ2}\\ 
   	\frac{dV}{ds} &=&  \alpha U^*(s)   \ -  y V(s)    \ -  \mu_{v} V(s) \label{sec2equ3}
   \end{eqnarray} 
 
 where $$x=\bigg(d_{1}  + d_{2}  + d_{3} +  d_{4} + d_{5}+ d_{6}\bigg)$$
 $$y=\bigg( b_{1}  + b_{2}  + b_{3}  +  b_{4}+ b_{5} + b_{6}\bigg)$$
 
 The viral load $V(s)$ at within-host scale provide the link between the within-host scale and the between-host scale. We use the within-host scale sub-model given by the model system $(2.14)-(2.16)$ to estimate the individual host infectiousness. To derive an expression for the area under the viral load curve we use the within-host sub-model  $(2.14)-(2.16)$. We denote the area under the viral load curve by $N_h$. The quantity $N_h$ measures the total amount of COVID-19 virus produced by an infected human throughout the period of host infectivity, and thus can be considered a proxy for individual host infectivity \cite{hadjichrysanthou2016understanding}. Let $d_1$ and $d_2$ be the times at which the viral load takes on the value of the detection limit at the beginning and end of infection. Then $d_2-d_1$ can be considered as the duration of host infectivity. Using the method described in \cite{garira2020development} the area under the viral load curve is given by, 
 \begin{equation}
     N_h = \frac{\alpha \int_{d_1}^{d_2}U^* ds}{y+\mu_v}
 \end{equation}
Where, $U^*$ is the infected cells. For the chosen parameter values $N_h$ will have a fixed numerical value.

 Now the within-host scale sub-model $(2.14)-(2.16)$ is reduced to a composed parameter
$N_h$ given by Equation $(2.17)$. This $N_h$ replaces $V(s)$ in the
between-host scale sub-model as follows: 
 \begin{eqnarray}
 	\frac{dS}{dt}& =&  \Lambda \ - \beta N_h S(t)I(t)  - \mu S(t)  \label{sec2equ1} \\
   	\frac{dE}{dt} &=& \beta N_h S(t)I(t) \ - (\mu+\pi+\gamma_1)E(t)    \label{sec2equ2}\\ 
   	\frac{dI}{dt} &=& \pi E(t)-(\mu+\gamma_2)I(t)-d N_h I(t)\label{sec2equ3} 
   \end{eqnarray}
 
 The reduced simplified model $(2.18)-(2.20)$ is now at a single time scale $t$ and is much easier to analyse now. We will study the dynamics of  COVID-19 disease using the reduced simplified model $(2.18)-(2.20)$ and also study the influence of within-host sub-model parameters on the spread of infection. 
 
 \subsection{\textbf{Well-Posedness of the Model}}
 The existence, the positivity, and the boundedness of the solutions of the proposed model $(2.18)-(2.20)$ need to be proved to ensure that the model has a mathematical and biological meaning. 

\subsubsection{Positivity and Boundedness}
\underline{\textbf{Positivity}}\textbf{:}

\begin{lem}
 Let $t_0 > 0$ and $S(t_0) > 0$, $E(t_0) > 0$, $I(t_0) > 0$ then the solution $S(t), E(t)$ and $I(t)$ of the system $(2.18)-(2.20)$ are positive for all $t\geq 0$.
 
proof:  From equation $(2.18)$ we have 

 \begin{equation*}
\begin{split}
\frac{dS}{dt}&= \Lambda - \beta N_h S I-\mu S  \\[4pt]
\frac{dS}{dt}& \geq -(\beta N_h S I -\mu) S \\
\frac{dS}{S}& =  (\beta N_h  I -\mu) dt\\
 \end{split}
\end{equation*}
 Integrating both sides from $t_0$ to $t$ we get
 
 $$S(t)\geq S(t_0) e^{-\int_{t_0}^{t} (\beta N_h I + \mu)}$$
  Therefore $S(t) > 0$ for all $t>0$.
  
  Also from equation $2.19$ we have,
 \begin{equation*}
\begin{split}
\frac{dE}{dt}&=  \beta N_h S I -(\mu+\pi + \gamma_1 )E  \\[4pt]
\frac{dE}{dt}& = \frac{\beta N_h S I E}{E} -(\mu+\pi + \gamma_1 )E  \\
\frac{dE}{E}& = \bigg( \frac{\beta N_h S I }{E} -(\mu+\pi + \gamma_1 )\bigg)dt\\
\frac{dE}{E}& = f(S, E,I) dt\\
 \end{split}
\end{equation*}

where  $$f(S, E,I)= \frac{\beta N_h S I }{E} -(\mu+\pi + \gamma_1 )$$

 Integrating both sides from $t_0$ to $t$ we get
 $$E(t)= E(t_0) e^{\int_{t_0}^{t} f(S,E,I)}$$
  Therefore $E(t) > 0$ for all $t>0$.

From the last equation $2.20$ we have

 \begin{equation*}
\begin{split}
\frac{dI}{dt}&= \pi E -(\mu+\gamma_2 + d N_h) I  \\[4pt]
\frac{dI}{dt}& \geq -(\mu+\gamma_2 + d N_h) I \\
\frac{dI}{I}& =  (\mu+\gamma_2 + d N_h)  dt\\
 \end{split}
\end{equation*}
 Integrating both sides from $t_0$ to $t$ we get
 
 $$I(t)\geq I(t_0) e^{-\int_{t_0}^{t} (\mu+\gamma_2+d N_h)}$$
  Therefore $I(t) > 0$ for all $t>0$.
  Hence, we conclude that all the solutions of the the system $(2.18)-(2.20)$ remain positive for any time $t>0$  provided that the initial conditions are positive. This establishes the positivity of the solutions of the system $(2.18)-(2.20)$.
\end{lem}

\underline{\textbf{Boundedness}}\textbf{:}\\
Let  $N(t) = S(t)+E(t)+I(t) $ \\
Now,  
\begin{equation*}
\begin{split}
\frac{dN}{dt} & = \frac{dS}{dt} +  \frac{dE}{dt}+  \frac{dI}{dt}  \\[4pt]
& = \Lambda -\mu(S(t) + E(t) + I(t)) -\gamma_1 E-\gamma_2 I - d N_h \\
& \leq \Lambda - \mu N(t)
\end{split}
\end{equation*}
Therefore,
$$\frac{dN}{dt} + \mu N(t) \leq \Lambda$$
 The integrating factor is given by $e^{\mu t}.$ Therefore, after integration we get,
$$N(t)\le \frac{\Lambda}{\mu}$$
Thus we have shown that the solutions of the system $(2.18)-(2.20)$ is bounded. 

Therefore, the biologically feasible region is given by the following set, 
\begin{equation*}
\Omega = \bigg\{\bigg(S(t), E(t), I(t)\bigg) \in \mathbb{R}^{3}_{+} : S(t)  + E(t) + I(t) \leq \frac{\Lambda}{\mu}, \ t \geq 0 \bigg\}
\end{equation*}

We summarize the above discussion on  boundedness by the following lemma. 

\begin{lem}
    The set
    $$\Omega = \bigg\{\bigg(S(t), E(t), I(t)\bigg) \in \mathbb{R}^{3}_{+} : S(t)  + E(t) + I(t) \leq \frac{\Lambda}{\mu}, \ t \geq 0 \bigg\}$$
    is a positive invariant and an attracting set for system $(2.18)-(2.20)$.
\end{lem}

\subsubsection{Existence and Uniqueness of Solution}

For the general first order ODE of the form
\begin{equation}
    \dot x=f(t,x) , \hspace{2cm}x(t_0)=x_0
    \end{equation}
One would have interest in knowing the answers to the following questions: \\
(i) Under what conditions  solution exists for the $(2.21)$?\\
(ii) Under what conditions unique solution exists for the system $(2.21)$?

We use the following theorem to established the existence and uniqueness of solution for our SEI  model $(2.18)-(2.20)$.

\begin{thm}
Let D denote the domain: \\
$$|t-t_0| \leq a, ||x-x_0|| \leq b, x=(x_1, x_2,..., x_n), x_0=(x_{10},..,x_{n0})$$
and suppose that $f(t, x)$ satisfies the Lipschitz condition:
\begin{equation}
||f(t,x_2)-f(t, x_1)|| \leq k ||x_2-x_1||
\end{equation}
and whenever the pairs $(t,x_1)$ and $(t, x_2)$ belong to the domain $D$ , where $k$ is used to represent a positive constant. Then, there exist a constant $\delta > 0$ such that there exists a unique (exactly one) continuous vector solution $x(t)$ of the
system $(2.21)$ in the interval $|t-t_0|\leq \delta $. It is important to note that condition $(2.22)$ is satisfied by requirement that:
$$\frac{\partial f_i}{\partial x_j}, i,j=1,2,.., n$$
be continuous and bounded in the domain D.
\end{thm}

\begin{thm}{Existence of Solution}\\
Let $D$ be the domain defined in above such that $(2.22)$ hold. Then there exist a solution of model system of equations $(2.18)-(2.20)$ which is bounded in the domain $D$.\\
proof\\
Let:
 \begin{eqnarray}
 	f_1& =&  \Lambda \ - \beta N_h S(t)I(t)  - \mu S(t)  \label{sec2equ1} \\
   	f_2 &=& \beta N_h S(t)I(t) \ - (\mu+\pi+\gamma_1)E(t)    \label{sec2equ2}\\ 
   f_3 &=& \pi E(t)-(\mu+\gamma_2)I(t)-d N_h I(t)\label{sec2equ3} 
   \end{eqnarray}
  We will show that 
  $$\frac{\partial f_i}{\partial x_j}, i,j=1,2,.., n$$
 is continuous and bounded in the domain D.
 
 From equation $(2.23)$ we have
 \begin{eqnarray*}
    \frac{\partial f_1}{\partial S}&=&-\beta N_h I-\mu, |\frac{\partial f_1}{\partial S}|=|-\beta N_h I-\mu| < \infty\\
 \frac{\partial f_1}{\partial E}&=0&, |\frac{\partial f_1}{\partial E}| < \infty\\
 \frac{\partial f_1}{\partial I}&=&-\beta N_h S, |\frac{\partial f_1}{\partial I}|=|-\beta N_h S| < \infty
  \end{eqnarray*}
  
   Similarly from equation $(2.24)$ we have
 \begin{eqnarray*}
    \frac{\partial f_2}{\partial S}&=&\beta N_h I, |\frac{\partial f_2}{\partial S}|=|\beta N_h I| < \infty\\
 \frac{\partial f_2}{\partial E}&=&-\mu-\pi-\gamma_1, |\frac{\partial f_2}{\partial E}| = |-(\mu+\pi+\gamma_1)|< \infty\\
 \frac{\partial f_2}{\partial I}&=&\beta N_h S, |\frac{\partial f_1}{\partial I}|=|\beta N_h S| < \infty
  \end{eqnarray*}
  
Finally from $(2.25)$ we have
  \begin{eqnarray*}
    \frac{\partial f_3}{\partial S}&=&0, |\frac{\partial f_3}{\partial S}| < \infty\\
 \frac{\partial f_3}{\partial E}&=&\pi, |\frac{\partial f_3}{\partial E}| = |\pi|< \infty\\
 \frac{\partial f_3}{\partial I}&=&-(\mu + \gamma_2+ d N_h), |\frac{\partial f_3}{\partial I}|=|-(\mu + \gamma_2+ d N_h)| < \infty
  \end{eqnarray*}
  
  Hence we have shown that all the partial derivatives are continuous and bounded. Therefore, Lipschitz condition $(2.22)$ is satisfied. Hence, by theorem 2.1 there exists a unique solution of system $(2.18)-(2.20)$ in the region $D$. 
\end{thm}

 \subsection{{\textbf{Stability Analysis}}}\vspace{.25cm}

 The basic reproduction number denoted by $R_0$ that gives the average number of secondary cases per primary case is calculated using the next generation matrix method \cite{diekmann2010construction} and the expression for $R_{0}$ for the system $(2.18)-(2.20)$ is given by \vspace{.5cm}\\

\begin{equation}
\mathbf{ R_{0}}= \mathbf{\frac{\beta N_h \pi \Lambda}{\mu (\mu+ \pi +\gamma_1) (\mu + \gamma_2 + d N_h)}} \label{sec3equ1}\\
\end{equation}

System  $(2.18)-(2.20)$ admits two equilibria  namely, the infection free equilibrium $E_{0}=\bigg(\frac{\omega}{\mu},0,0 \bigg)$ and the infected equilibrium $E_{1}=(S^*, E^*, I^*)$  where,

$$\hspace{-1cm}S^*=\frac{\Lambda}{(\beta N_h I^* + \mu)}$$

$$E^*=\frac{(\mu + \gamma_2 + d N_h)I^*}{\pi}$$

$$I^*=\frac{\mu (R_0 - 1)}{\beta N_h}$$

Since negative population does not make sense, the existence condition for the infected equilibrium point $E_1$ is that $R_0 > 1$.

\subsubsection{\textbf{Stability Analysis of $E_0$}}
 We analyse the stability of  equilibrium points $E_0$. This is done based on the nature of the eigenvalues of the jacobian matrix evaluated at $E_0$.

The jacobian matrix of the system $(2.18)-(2.20)$ at the infection free equilibrium $E_0$ is given by, \\

\begin{equation*}
J_{E_{0}} = 
\begin{pmatrix}
-\mu & 0 & \frac{-\beta N_h \Lambda}{\mu} \\
0 & -(\mu + \pi +\gamma_1) & \frac{\beta N_h \Lambda }{\mu} \\
0 & \pi & -(\mu + \gamma_2 + d N_h)
\end{pmatrix}
\end{equation*} \\
 
 The characteristic equation of $J_{E_0}$ is given by,
 \begin{equation}
 (-\mu-\lambda)\bigg(\lambda^2+(2\mu+ \pi + \gamma_1+\gamma_2 + d N_h)\lambda -(R_0 -1) (\mu + \pi + \gamma_1)(\mu + \gamma_2 + d N_h)\bigg) \label{111b}
 \end{equation}
 One of the eigenvalue of characteristic equation $(\ref{111b})$ is $-\mu$ and the other two are the roots of the following quadratic equation.\\
 \begin{equation}
     \lambda^2+(2\mu+ \pi + \gamma_1+\gamma_2 + d N_h)\lambda -(R_0 -1) (\mu + \pi + \gamma_1)(\mu + \gamma_2 + d N_h)
 \end{equation}
We see that when $R_0 < 1$ both the eigenvalues of equation $(2.28)$ is negative. Therefore, $E_0$ remains locally asymptotically stable whenever $R_0 < 1$. When $R_0 > 1$ the quadratic equation $(2.28)$ has one positive and one negative root. Therefore the characteristic equation $(2.27)$ has two negative and one positive root. Hence $E_0$ becomes unstable in this case. We summarize the local asymptotic stability of infection free equilibrium $E_0$ in the following theorem.

\begin{thm}
      The infection free equilibrium point $E_0$ of  system $(2.18)-(2.20)$ is locally asymptotically stable provided $R_0 < 1$. If $R_0$ crosses unity $E_0$ loses its stability and becomes unstable.   
\end{thm}

 \newpage
{\flushleft{  \textbf{Global Stability of  $E_{0}$ }}}

To establish the global stability of the infection free equilibrium $E_{0}$ we make use of the method discussed in Castillo-Chavez \textit{et al}  {\cite{GLB}}.\\ 

\begin{thm}

Consider the following general system,\\
\begin{equation}
\begin{split}
\frac{dX}{dt} & = F(X,Y) \\[6pt] \label{sec3equ5}
\frac{dY}{dt} & = G(X,Y) \qquad   
\end{split}
\end{equation}
where $ X $ denotes the uninfected population compartments and $Y$ denotes the infected population compartments including latent, infectious etc. Here the function $G$ is such that it satisfies $G(X,0)=0$. Let $U_{0} = (X_{0}, \bar{0}) $ denote the equilibrium point of the above general system.

If the following  two conditions are satisfied then the infection free equilibrium point $U_{0}$ is globally asymptotically stable for the above general system provided $R_{0}<1$

$A_{1}$: For the subsystem $\frac{dX}{dt} = F(X,0)$, $X_{0}$ is globally asymptotically stable.

$A_{2}$: The function $G=G(X,Y)$ can be written as $G(X,Y)=AY - \widehat{G}(X,Y)$, where $\widehat{G_{j}}(X,Y) \geq 0$ $\forall \; (X,Y)$ in the biologically feasible region $\Omega$ for j=1,2 and $A = D_{Y}G(X,Y)$ at $(X_{0},\bar{0})$ is a M-matrix(matrix with non-negative off diagonal element).
\end{thm}
\vspace{.5cm}
\begin{thm}
The infection free equilibrium point $E_0$ of the system $(2.18)-(2.20)$ is globally asymptotically stable whenever $R_0 < 1$\\
proof:\\
We will prove global stability of $E_{0}=(\frac{\omega}{\mu},0,0)$ of system $(2.18)-(2.20)$ by showing that system  $(2.18)-(2.20)$ can be written as the above general form and both the conditions $A_{1}$ and $A_{2}$ are satisfied .\\
Comparing the above general system (\ref{sec3equ5}) to the system $(2.18)-(2.20)$ the functions $F$ and $G$ are given by
\\
$$   F(X,Y) = \Lambda - \beta N_h S I - \mu S$$
$$\hspace{1.3cm}G(X,Y) = \bigg( \beta N_h S I-(mu+\pi + \gamma_1) E, \; \pi E-(\mu + \gamma_2 +d N_h) I\bigg)$$ 

\noindent
\\ 
where $X=S$ and $Y=(E, I)$ \\
The disease free equilibrium point is $U_{0} = (X_{0}, \bar{0}), $ where,
\\
$$X_{0} = \frac{\Lambda}{\mu} \quad \text{and} \quad \bar{0} = (0, 0)$$

\noindent
\\
From the stability analysis of $E_{0}$, we know that  $U_{0}$ is locally asymptotically stable iff $R_{0} < 1$. Clearly, we see that $G(X,\bar{0}) = (0,\bar{0})$. Now, we show that $ X_{0} = (\frac{\Lambda}{\mu})$ is globally asymptotically stable for the subsystem 

\begin{equation}
    \frac{dS}{dt}=F(S,\bar{0})=\omega - \mu S  \label{sec3equ6}
\end{equation}
\noindent
The integrating factor is $e^{\mu t}$ and therefore after performing integration on the above equation (\ref{sec3equ6}) we get, 
$$S(t)e^{\mu t}=\frac{\Lambda e^{\mu t}}{\mu }  + c$$
\\\noindent
As $t \rightarrow \infty$ we get, 
$$S(t)= \frac{\omega}{\mu}$$
\noindent
which is independent of c. This independency implies that $X_{0}=\frac{\Lambda}{\mu_{1}}$ is globally asymptotically stable for the subsystem $\frac{dS}{dt}=\Lambda - \mu S$. So, the assumption  $A_{1}$ is satisfied.\\

\noindent
Now, we will show that assumption $A_{2}$  holds. First, we will find the matrix $A$. As per the theorem, $A = D_{Y}G(X,Y)$ at $X=X_{0}$ and $Y=\bar{0}$. Now\\ 
\vspace{2mm}
\[ D_{Y}G(X,Y) = 
\begin{bmatrix}
-(\mu + \pi + \gamma_1) & \beta N_h S \\[6pt]
\pi & -(\mu + \gamma_2 + d N_h)
\end{bmatrix}\]

\noindent
\\ 
At $X=X_{0}$ and $Y=\bar{0}$, we obtain,
\\ \vspace{1mm}
\[A = 
\begin{bmatrix}
-(\mu + \pi + \gamma_1) & \frac{\beta N_h \Lambda}{\mu} \\[6pt]
\pi & -(\mu + \gamma_2 + d N_h)
\end{bmatrix}\]
\vspace{2mm}
\noindent
\\ 
Clearly, matrix $A$ has non-negative off-diagonal elements. Hence, $A$ is a M-matrix. Using 
$\widehat{G}(X,Y) = AY - G(X,Y)$, we get,
\vspace{2mm}
\[\widehat{G}(X,Y) \;\; = \;\;
\begin{bmatrix}
\widehat{G_{1}}(X,Y) \\[6pt]    
\widehat{G_{2}}(X,Y)    
\end{bmatrix} 
\;\; = \;\;
\begin{bmatrix}
\beta N_h I(\frac{\Lambda}{\mu}-S) \\[6pt]
0
\end{bmatrix} \] \\

Hence $\widehat{G_{1}}(X,Y) = \beta N_h I(\frac{\Lambda}{\mu}-S) \geq 0 $ because $S(t)\le\frac{\Lambda}{\mu}$ and  $\widehat{G_{2}}(X,Y) =0$ \\
Thus both the assumptions $A_{1}$ and $A_{2}$ are satisfied and therefore infection free equilibrium point $E_{0}$ is globally asymptotically stable provided $R_{0} < 1.$ 

\end{thm}

\subsubsection{Stability Analysis of $E_1$}
The jacobian matrix of the system $(2.18)-(2.20)$ at $E_1$ is given by, \\

\begin{equation*}
J = 
\begin{pmatrix}
-(\beta N_h I^*+\mu) & 0 & -\beta N_h S \\
\beta N_h I^* & -(\mu+\pi+\gamma_1) & \beta N_h S^* \\
0 & \pi & -(\mu+\gamma_2 + d N_h)
\end{pmatrix}
\end{equation*} \\

The characteristic equation of the jacobian $J$ evaluated at $E_{1}$ is given by, \\
\begin{equation}
\lambda^3 + A_1 \lambda^2 + B_1\lambda + C_1 = 0 \label{sec3equ4}  \\
\end{equation}
where 
$$A_1=3\mu + \pi + \beta N_h I^* + \gamma_1+\gamma_2 + dN_h$$
$$B_1=(\mu+\beta N_h I^*)(2\mu + \pi + \gamma_1+\gamma_2 + d N_h)+(\mu + \pi +\gamma_1)(\mu + \gamma_2+d N_h)-\beta N_h \pi S^*$$
$$C_1=(\mu + \beta N_h I^*)\bigg((\mu + \pi +\gamma_1)(\mu + \gamma_2+d N_h)-\beta N_h \pi S^*\bigg)+ \beta^2 N_h^2 S^* I^* \pi$$

Clearly, $A_1 > 0$. By Routh- Hurwitz criterion, all the roots of characteristic equation $(2.31)$ are negative iff $C_1 > 0$ and $A_1 B_1-C_1 > 0$.\\
Simplifying the expression for $C_1$ we get,
$$C_1=\mu(\mu + \pi + \gamma_1)(\mu + \gamma_2 + d N_h)(R_0-1)$$
Therefore, $C_1 > 0$ iff $R_0 > 1$.
Hence we conclude that the infected equilibrium point $E_{1}$ exists and remains locally asymptotically stable provided $R_{0} > 1.$ and $(A_1 B_1-C_1) >0$. In the following theorem we summarize the above discussion on the stability of $E_1$.

\begin{thm}
There exists a unique infected equilibrium point $E_1$ of the system $(2.18)-(2.20)$ if the following conditions are satisfied:\\
(i) $R_0 > 1$ \\
(ii) $(A_1 B_1-C_1)>0$
\end{thm}

\subsection{\textbf{Bifurcation Analysis}} 

 We now use the method given by Chavez and Song in  \cite{BIF} to do the bifurcation  analysis.

\vspace{.25cm}

\begin{thm}

	Consider a system, 
	$$\frac{dX}{dt}=f(X,\phi)$$
	where $X \in \mathbb{R}^{n}$, $\phi \in \mathbb{R}$ is the bifurcation parameter and $f : \mathbb{R}^{n} \times \mathbb{R} \rightarrow \mathbb{R}^{n} $ where $f \in \mathbb{C}^2 (\mathbb{R}^n, \mathbb{R})$. Let $\bar{0}$ be the equilibrium point of the system such that $f(\bar{0},\phi) = \bar{0}, \forall \; \phi \in \mathbb{R}$. 
	Let the following conditions hold :
	
	\begin{enumerate}
		\item For the matrix $A = D_{X}f(\bar{0},0)$, zero is the simple eigenvalue and all other eigenvalues have negative real parts.
		
		\item Corresponding to zero eigenvalue, matrix A has non-negative right eigenvector, denoted as $u$ and non-negative left eigenvectors, denoted as $v$.
	\end{enumerate}
	
	\noindent
	Let $f_{k}$ be the $k^{th}$ component of $f$. Let $a$ and $b$ be defined as follows -
	
	$$ a = \sum_{k,i,j=1}^{n} \Bigg[ v_{k} w_{i} w_{j} \bigg(\frac{\partial^{2} f_{k}}{\partial x_{i} \partial x{j}} (\bar{0},0)\bigg) \Bigg]$$
	
	$$ \hspace*{-0.65cm} b = \sum_{k,i=1}^{n} \Bigg[ v_{k} w_{i} \bigg(\frac{\partial^{2} f_{k}}{\partial x_{i} \partial \phi}(\bar{0},0)\bigg) \Bigg]$$

	\noindent
	Then local dynamics of the system near the equilibrium point $\bar{0}$ is totally determined by the signs of $a$ and $b$. Here are the following conclusions :
	
	\begin{enumerate}%[\label=(\roman*)]
		\item If $a > 0$ and $b > 0$, then whenever $\phi < 0$ with $\mid \phi \mid \ll 1$, the equilibrium $\bar{0}$ is locally asymptotically stable, and moreover there exists a positive unstable equilibrium. However when $0 < \phi \ll 1$, $\bar{0}$ is an unstable equilibrium and there exists a negative and locally asymptotically stable equilibrium.
		
		\item If $ a < 0$, $b < 0$, then whenever $\phi < 0$ with $\mid \phi \ll 1 $, $\bar{0}$ is an unstable equilibrium whereas if $0 < \phi \ll 1$, $\bar{0}$ is locally asymptotically stable equilibrium and there exists a positive unstable equilibrium.
		
		\item If $a > 0$, $b < 0$, then whenever $\phi < 0$ with $\mid \phi \mid \ll 1$, $\bar{0}$ is an unstable equilibrium, and there exists a locally asymptotically stable negative equilibrium. However if $0 < \phi \ll 1$, $\bar{0}$ is stable, and a there appears a positive unstable equilibrium.
		
		\item If $a < 0$, $b > 0$, then whenever $\phi$ changes its value from negative to positive, the equilibrium $\bar{0}$ changes its stability from stable to unstable. Correspondingly a negative  equilibrium, unstable in nature, becomes positive and locally asymptotically stable.
	\end{enumerate}

\end{thm}
\underline{\textbf{Applying the Theorem 2.7 to our system $(2.18)-(2.20)$ }}\textbf{: }
\noindent

In our case, we have $x = (S, E, I) \in \mathbb{R}^3$ where $x_{1} = S$, $x_{2} = E$ and $x_{3} = I$. Let us consider $\beta$ (transmission rate of the infection) to be the bifurcation parameter.\\
We know that, 
$$R_{0}=\frac{\beta N_h \Lambda \pi }{\mu(\mu+\pi+\gamma_1)(\mu +\gamma_2+d N_h)}$$
Therefore we have, \\

\begin{equation*}
\beta = \frac{ R_{0}\mu(\mu+\pi+\gamma_1)(\mu +\gamma_2+d N_h)}{\beta N_h \Lambda \pi}
\end{equation*}
\noindent
\\
Let $\beta = \beta^{*}$ at $R_{0}=1$. So, we have,
\vspace{2mm}
\begin{equation*}
\beta^* = \frac{\mu(\mu+\pi+\gamma_1)(\mu +\gamma_2+d N_h)}{\beta N_h \Lambda \pi}
\end{equation*}

\noindent
With $x=(x_{1},x_{2},x_3)=(S,I,V)$
system $(2.18)-(2.20)$ can be written as follows :

\begin{align*}
\frac{dx_{1}}{dt} &= \Lambda-\beta N_h x_{1}x_{3}-\mu x_{1} = f_{1}\\[6pt] 
\frac{dx_{2}}{dt} &=  \beta N_h x_{1}x_{3}-(\mu+\pi +\gamma_1)x_{2} = f_{2}\\[6pt]
\frac{dx_{3}}{dt} &= \pi x_{2}-(\mu+\gamma_2 + d N_h) x_{3} \quad= f_{3}
\end{align*}

\noindent
The disease free equilibrium point $E_{0}$ is given by,

$$x^* = \bigg(\frac{\Lambda}{\mu}, 0, 0\bigg) = (x_{1}^*, x_{2}^*, x_{3}^*)$$
\vspace{2mm}
\noindent
\\ 
Clearly, $f(x^*,\beta) = 0, \; \forall \; \beta \in \mathbb{R}$, where $ f = (f_{1},f_{2}, f_{3})$. Let $D_{x}f(x^*,\beta^*)$ denote the Jacobian matrix of the above system at the equilibrium point $x^*$ and $R_{0} = 1$. Now we see that,
\vspace{2mm}

\[
D_{x}f(x^*,\beta^*) =
\begin{bmatrix}
-\mu & 0 & \frac{-\beta^* N_h\Lambda}{\mu} \\[6pt]
0 & -(\mu+\pi +\gamma_1) & \frac{\beta^* N_h \Lambda}{\mu} \\[6pt]
0 & \pi & -(\mu + \gamma_2 + d N_h)
\end{bmatrix}
\]
\vspace{2mm}
\noindent
\\
The characteristic polynomial of the above matrix $D_{x}f(x^*,\beta^*)$ is given by, 

\begin{equation}
    (-\mu-\lambda)\bigg[((\mu + \pi + \gamma_1)+\lambda)((\mu + \gamma_2 + d N_h)+\lambda) - \bigg(\frac{\alpha \beta^* N_h \Lambda \pi}{\mu}\bigg)\bigg] = 0  \label{sec3equ8}
 \end{equation}

\vspace{2mm}
\noindent
Hence, we obtain the first eigenvalue  of (\ref{sec3equ8})
as $$\boldsymbol{\lambda_{1}= -\mu < 0} $$
\noindent
The other eigenvalues $\lambda_{2,3}$ of (\ref{sec3equ8}) are the solutions of the following equation,\\ 
\begin{equation}
    \lambda^2 +\bigg(2\mu + \pi + \gamma_1 +\gamma_2 + d N_h\bigg)\lambda+(\mu + \pi + \gamma_1)(\mu + \gamma_2 + d N_h)-\frac{\beta^* N_h \pi \Lambda}{\mu}=0  \label{sec3equ9}
\end{equation}

substituting the expression for  $\beta^*$ in (\ref{sec3equ9}) we get, \\
\begin{equation}
    \lambda^2 +\bigg(2\mu+ \pi + \gamma_1 +\gamma_2 + d N_h\bigg)\lambda = 0 \label{sec3equ10}
\end{equation}

The  eigen values of (\ref{sec3equ10}) are $\lambda_{2}=0$ and $\lambda_{3} = -(2\mu+ \pi + \gamma_1 +\gamma_2 + d N_h)$

\noindent
Hence, the matrix $D_{x}f(x^*,\beta^*)$ has zero as its simple eigenvalue and all other eigenvalues with negative real parts. Thus, the condition 1 of the theorem 2.7 is satisfied.
\vspace{1cm}\\
Next, for proving condition 2, we need to find the right and left eigenvectors of the zero eigenvalue ($\lambda_{2}$). Let us denote the right and left eigen vectors by $\boldsymbol{w}$ and $\boldsymbol{v}$ respectively. To find $w$, we use $(D_{x}f(x^*,\beta^*) - \lambda_{2} I_{d})w = 0 $, which implies that 
\\

\[\
\begin{bmatrix}
-\mu & 0 & \frac{-\beta^* N_h \Lambda}{\mu} \\[6pt]
0 & -(\mu+\pi +\gamma_1) & \frac{\beta^* N_h \Lambda}{\mu} \\[6pt]
0 & \pi & -(\mu + \gamma_2 + d N_h)

\end{bmatrix}
\begin{bmatrix}
w_{1} \\[6pt]
w_{2} \\[6pt]
w_{3} \\[6pt]
\end{bmatrix}
\;\; = \;\;
\begin{bmatrix}
0 \\[6pt]
0 \\[6pt]
0 
\end{bmatrix}
\]

\vspace{2mm}
\noindent
\\ 
where $ w = (w_{1}, w_{2}, w_{3})^T$. As a result, we obtain the system of simultaneous equations as follows :

\begin{equation}
    -\mu w_{1} - \frac{\beta^* N_h \Lambda}{\mu}w_{3} = 0 \label{sec3equ11}
\end{equation}
 
 \begin{equation}
 -(\mu+\pi +\gamma_1)w_{2} + \frac{\beta^* N_h\Lambda}{\mu}w_{3} = 0 \label{sec3equ12}
\end{equation}

 \begin{equation}
 \pi w_{2} - (\mu + \gamma_2 + d N_h) w_{3} = 0 \label{sec3equ13}
\end{equation}
\noindent
\\
By choosing $w_{3} = \mu$ in the above simultaneous equation  (\ref{sec3equ11})-(\ref{sec3equ13}) we obtain 
$$ w_{2} = \frac{\beta^* N_h \Lambda}{(\mu+ \pi+\gamma_1)} \;\;\; \text{and} \;\;\;{w_{1} = -\frac{\beta^* N_h \Lambda}{\mu}} $$

\noindent
Therefore, the right eigen vector of zero eigenvalue is given by
\begin{equation*}
\boldsymbol{w = \bigg(-\frac{\beta^* N_h\Lambda}{\mu}, \; \frac{\beta^* N_h\Lambda}{(\mu+\pi +\gamma_1)}, \; \mu \bigg)} 
\end{equation*}

\noindent

Similarly, to find the left eigenvector $v$, we use $v(D_{x}f(x^*,\beta^*) - \lambda_{2} I_{d}) = 0 $, which implies that

\[
\begin{bmatrix}
v_{1} & v_{2} & v_{3}
\end{bmatrix}
\begin{bmatrix}
-\mu & 0 & \frac{-\beta^* N_h \Lambda}{\mu} \\[6pt]
0 & -(\mu+\pi +\gamma_1) & \frac{\beta^* N_h \Lambda}{\mu} \\[6pt]
0 & \pi & -(\mu + \gamma_2 + d N_h)

\end{bmatrix}
\;\; = \;\;
\begin{bmatrix}
0 & 0 & 0
\end{bmatrix}
\]

\noindent

where $v = (v_{1}, v_{2}, v_{3})$. The simultaneous equations obtained thereby are as follows :

\begin{equation}
     -\mu v_{1} = 0  \label{sec3equ14}
\end{equation}

\begin{equation}
 -(\mu + \pi + \gamma_1)v_{2} + \pi v_{3} = 0    \label{sec3equ15}
\end{equation}

\begin{equation}
 \frac{-\beta^* N_h\Lambda}{\mu}v_{1}+ \frac{\beta^* N_h \Lambda}{\mu}v_{2} -(\mu + \gamma_2 + d N_h)v_{3} = 0   \label{sec3equ16} 
\end{equation}
\vspace{2mm}
\noindent
\\
Therefore solving the above simultaneous equation (\ref{sec3equ14} - \ref{sec3equ16}) we obtain $v_{1} = 0$.\\

By choosing $v_{2}=1$ we get $$ v_{3}= \frac{\beta^* N_h \Lambda}{(\mu(\mu + \gamma_2 + d N_h))} $$

\noindent
\\
Hence, the left eigen vector is given by 

\begin{equation*}
\boldsymbol{v = \bigg(0, \; 1, \; \frac{\beta^* N_h \Lambda}{(\mu(\mu + \gamma_2 + d N_h))}\;\bigg)} 
\end{equation*}

\noindent
Now, we need to find $a$ and $b$. As per the Theorem 2.7, $a$ and $b$ are given by
\\ \vspace{2mm}
$$ \hspace*{-10mm} \boldsymbol{a = \sum_{k,i,j=1}^{3} \Bigg[v_{k}u_{i}u_{j} \bigg(\frac{\partial^2 f_{k}}{\partial x_{i} \partial x_{j}}(x^*, \beta^*)\bigg) \Bigg]}$$
\vspace{2mm}
$$ \hspace*{-15mm} \boldsymbol{b =\sum_{k,i =1}^{3} \Bigg[v_{k}u_{i}\bigg(\frac{\partial^2 f_{k}}{\partial x_{i} \partial \beta}(x^*, \beta^*)\bigg)\Bigg]}$$

\noindent
\\
Expanding the summation in the expression for $a$, it reduces to

$$ a = w_{1}w_{3} \frac{\partial^2 f_{2}}{\partial x_{1} \partial x_{3}} \thinspace + \thinspace \thinspace  w_{3}w_{1}\frac{\partial^2 f_{2}}{\partial x_{3} \partial x_{1}} $$

\noindent
\\ where partial derivatives are found at $(x^*, \beta^*)$. Now
$$ \frac{\partial^2 f_{2}}{\partial x_{1} \partial x_{3}}(x^*, \beta^*) = \beta^* N_h \qquad \frac{\partial^2 f_{2}}{\partial x_{3} \partial x_{1}}(x^*, \beta^*) = \beta^*N_h $$

\noindent
Substituting these partial derivatives along with $w$  in the expression of $a$, we get,
$$ \boldsymbol{a = -2\beta^{*2} N_h^2\Lambda < 0}$$

\noindent
Next, expanding the summation in the expression for $b$, we get,
$$b = v_{2}w_{2}\bigg(\frac{\partial^2 f_{2}}{\partial x_{1} \partial \beta}(x^*, \beta^*)\bigg) + v_{2}w_{3}\bigg(\frac{\partial^2 f_{2}}{\partial x_{3} \partial \beta}(x^*, \beta^*)\bigg)$$

\noindent
\\
Now
$$ \frac{\partial^2 f_{2}}{\partial x_{3} \partial \beta}(x^*, \beta^*) = \frac{N_h \lambda}{\mu} $$
$$ \frac{\partial^2 f_{2}}{\partial x_{1} \partial \beta}(x^*, \beta^*) = 0 $$
\noindent
substituting in the expression of $b$ we get,
$$ \boldsymbol{b = N_h\Lambda > 0}$$

\noindent
Hence $a < 0$ and $b>0$. We notice that condition (iv) of the theorem 2.7 is satisfied. Hence, we conclude that the system undergoes bifurcation at $\beta = \beta^*$ implying $R_{0} = 1$.

\noindent
Thus, we conclude that when $R_{0} < 1$, there exists a unique disease free equilibrium  which is globally asymptotically stable  and negative infected equilibrium which is unstable . Since negative values of population is not practical, therefore we ignore it in this case. Further, as $R_{0}$ crosses unity from below, the disease free equilibrium point loses its stable nature and become unstable, the bifurcation point being at $\beta = \beta^*$ implying $R_{0}=1$ and there appears a positive locally asymptotically stable infected equilibrium point. There is an exchange of stability between disease free equilibrium and infected equilibrium at $R_{0} = 1$. Hence, a forward bifurcation (\textbf{trans-critical bifurcation}) takes place at the break point $\beta = \beta^*$. We summarize the above discussion on bifurcation by the following theorem.

\begin{thm}
As $R_0$ crosses unity the disease free equilibrium changes its stability from
stable to unstable and there exists a locally asymptotically infected equilibrium when
$R_0 > 1$ i.e . direction of bifurcation is forward (transcritical) at $R_0 > 1$.

\end{thm}

 \subsection{\textbf{Sensitivity and Elasticity}}
 The Basic Reproduction number denoted by $R_0$ is one of the most important quantity in any infectious disease models. The expression for $R_0$ for the reduced multi-scale model $(2.18)-(2.20)$ is given by,
$$R_{0}=\frac{\beta N_h \Lambda \pi }{\mu(\mu+\pi+\gamma_1)(\mu +\gamma_2+d N_h)}$$

To determine best control measures, knowledge of the relative importance of the different factors responsible for
transmission is useful. Initially disease transmission is related to $R_0$ and sensitivity predicts which parameters have a high impact on $R_0$. The sensitivity index of $R_0$ with respect to a parameter $\mu$ is $\frac{\partial R_0}{\partial \mu}$ : Another measure is the elasticity index
(normalized sensitivity index) that measures the relative change of $R_0$ with respect to $\mu$, denoted by $\phi_{\mu}^{R_0}$, and defined as
 $$\phi_{\mu}^{R_0}=\frac{\partial R_0}{\partial \mu} \frac{\mu}{R_0}$$
 
 The sign of the elasticity index tells whether $R_0$ increases (positive sign) or decreases (negative sign) with the parameter;
whereas the magnitude determines the relative importance of the parameter \cite{berhe2019parameter, van2017reproduction}. These indices can guide control by indicating
the most important parameters to target, although feasibility and cost play a role in practical control strategy. If $R_0$ is known
explicitly, then the elasticity index for each parameter can be computed explicitly, and evaluated for a given set of parameters.

The elasticity index of $R_0$ with the model parameters is given by,\\
$$\phi_{\beta}^{R_0}=\frac{\partial R_0}{\partial \beta} \frac{\beta}{R_0}=1$$
 $$\phi_{\Lambda}^{R_0}=\frac{\partial R_0}{\partial \Lambda} \frac{\Lambda}{R_0}=1$$
 $$\phi_{\pi}^{R_0}=\frac{\partial R_0}{\partial \pi} \frac{\pi}{R_0}=1$$
 $$\phi_{\mu}^{R_0}=\frac{\partial R_0}{\partial \mu} \frac{\mu}{R_0}=\frac{-\beta N_h \Lambda \pi\bigg(3 \mu^2+ 2 \mu (\pi + \gamma_1 +\gamma_2 + d N_h) + (\gamma_2 +d)(\pi + \gamma_1)\bigg)}{\bigg(\mu^3 +\mu^2 (\pi + \gamma_1+\gamma_2 + d N_h)+ \mu (\gamma_2+d)(\pi + \gamma_1)\bigg)}\frac{\mu}{R_0}$$
 $$\phi_{\gamma_1}^{R_0}=\frac{\partial R_0}{\partial \gamma_1} \frac{\gamma_1}{R_0}=-\frac{\gamma_1}{(\mu +\pi + \gamma_1)}$$
  $$\phi_{\gamma_2}^{R_0}=\frac{\partial R_0}{\partial \gamma_2} \frac{\gamma_2}{R_0}=-\frac{\gamma_2}{(\mu  + \gamma_1+ d N_h)}$$
   $$\phi_{N_h}^{R_0}=\frac{\partial R_0}{\partial N_h} \frac{N_h}{R_0}=-\frac{\mu + \gamma_2}{(\mu+ \gamma_2  + d N_h)}$$
    $$\phi_{d}^{R_0}=\frac{\partial R_0}{\partial d} \frac{d}{R_0}=-\frac{d N_h}{(\mu+ \gamma_2  + d N_h)}$$

  The elasticity indices of the parameters of $R_0$ are given in Table 2. The elastic index of parameters $\beta, \pi, \Lambda$ and $N_h$ are positive and the remaining are negative. This implies that the increase in the values of these parameters increases $R_0$, whereas increase in the values of parameters $\mu, \gamma_1$ and $\gamma_2$ decreases $R_0$. For parameter $\beta$, $\phi_{\beta}^{R_0}=1$ implies an increase (decrease) of $\beta$ by $y \%$ increases (decreases) $R_0$ by the same percentage. From table 2 we see that the basic reproduction number is most sensitive to the the parameters $\beta, \pi$ and $\Lambda$. The implication of this is that an increase in the transmission rate increases the spread of the disease in the community.

  \begin{table}[ht!]
     	\caption{Elasticity Indices of $R_0$}
     	\centering % used for centering table
     	\begin{tabular}{|l|l|l|} % centered columns (4 columns)
     		\hline\hline
     		
     		\textbf{Parameters} &  \textbf{Elastic Index}& \textbf{Value} \\  % inserts table
     		%heading
     		\hline\hline % inserts single horizontal line
     	  $\beta$ & $\phi_{\beta}^{R_0} $&1 \\
     	\hline\hline
      $\pi$ & $\phi_{\pi}^{R_0}$ & 1\\
     	\hline\hline
     	 $\Lambda$ & $\phi_{\Lambda}^{R_0} $& 1\\
     	
     		\hline\hline
     		 $\mu$ & $\phi_{\mu}^{R_0} $& -0.2785\\
     		\hline\hline
     		 $\gamma_1$ & $\phi_{\gamma_1}^{R_0} $& -0.3196\\
     		 \hline\hline
     		 $\gamma_2$ & $\phi_{\gamma_2}^{R_0} $& -0.00082\\
     		 \hline\hline
     		 $N_h$ & $\phi_{N_h}^{R_0} $& 0.0018\\
     		 \hline\hline
     		 $d$ & $\phi_{d}^{R_0} $& -0.99\\
     		\hline\hline

     	\end{tabular}
     \end{table} \vspace{.25cm}

   \newpage
   
   \subsection{\textbf{Numerical Illustrations}}
 \subsubsection{\textbf{Numerical Illustrations of the stability of equilibrium points}}
 
Now we numerically illustrate the stability of the equilibrium points admitted by the system $(2.18)-(2.20)$. The simulation is done using matlab software and ode solver ode45 is used to solve the system of equations.  The parameter values of the within-host sub-model used in simulation is taken from \cite{chhetri2020within}. These within-host parameter values are used in calculating the area under the viral load curve $N_h$. The parameter values of the reduced multi-scale mode $(2.18)-(2.20)$ are taken from \cite{mwalili2020seir, samui2020mathematical}. All the parameter values are listed in table 4. These parameter values are used in illustrating the stability of the equilibrium points admitted by the system $(2.18)-(2.20)$. We also illustrate the influence of key within-host scale sub-model parameters on the between-host scale sub-model variables.

\begin{table}[ht!]
     	\caption{parameter Values}
     	\centering % used for centering table
     	\begin{tabular}{|l|l|l|} % centered columns (4 columns)
     		\hline\hline
     		
     		\textbf{Symbols} &  \textbf{Values}& \textbf{Source} \\  % inserts table
     		%heading
     		\hline\hline % inserts single horizontal line
     	  $\omega$ & $2$ & \cite{chhetri2020within} \\
     	\hline\hline
      $k$ & $0.05$ & \cite{chhetri2020within}\\
     	\hline\hline
     	 $\mu_c$ & $0.1 $& \cite{chhetri2020within}\\
     	\hline\hline
     		 $\mu_v$ & $0.1$ &\cite{chhetri2020within} \\
     		\hline\hline
     		 $\alpha$ & $0.24 $& \cite{mwalili2020seir} \\
     		 \hline\hline
     		 $d_1, d_2, d_3,d_4,d_5, d_6$ & $0.027,0.22,0.1,0.428,0.01,0.01$&  \cite{chhetri2020within} \\
     		 \hline\hline
     	 $b_1, b_2, b_3, b_4,b_5, b_6$ & $0.1,0.1,0.08,0.11,0.01,0.07$&  \cite{chhetri2020within} \\
     		\hline\hline
     		 $\Lambda$ & $\mu N(0) $& \cite{mwalili2020seir}\\
     		 \hline\hline
     		 $\beta$& 0.0115& \cite{samui2020mathematical}\\
     		  \hline\hline
     		 $\mu$ & 0.062 &\cite{mwalili2020seir}  \\
     		  \hline\hline
     		  $\pi$ & 0.09 & \cite{samui2020mathematical} \\
     		   \hline\hline
     		   $d$ & 0.0018 & \cite{samui2020mathematical}\\
     		    \hline\hline
     		    $\gamma_1, \gamma_2$ & 0.05, 0.0714 & \cite{samui2020mathematical}\\
     		     
     		 \hline\hline
     			
     	\end{tabular}
     \end{table} \vspace{.25cm}

The initial values used in the simulation is given in table 5.

  \begin{table}[ht!]
     	\caption{Initial Values of the Variables}
     	\centering % used for centering table
     	\begin{tabular}{|l|l|l|} % centered columns (4 columns)
     		\hline\hline
     		
     		\textbf{Variable} &  \textbf{Initial Values}& \textbf{Source} \\  % inserts table
     		%heading
     		\hline\hline % inserts single horizontal line
     	  $U(s)$ & $3.2*10^5$ & \cite{chhetri2020within} \\
     	\hline\hline
      $U^*(s)$ & $0$ & \cite{chhetri2020within}\\
     	\hline\hline
     	 $V(s)$ &$ 5.2 $& \cite{chhetri2020within}\\
     	    \hline\hline
     		 $S(t)$ & $1000 $& Assumed\\
     		\hline\hline
     		 $E(t)$ & $100 $& Assumed\\
     		 \hline\hline
     		 $I(t)$ & $50$& Assumed\\
     		 \hline\hline

     	\end{tabular}
     \end{table} \vspace{.25cm}

With the  parameter values from table 3, the area under the viral load curve $N_h$ is found to be $3.3759\times10^4$. The value of basic reproduction number $\boldsymbol{R_{0}}$ with $\beta=0.00115$,  $\mu=0.72$ and other parameters from table 3 is calculated to be  $\boldsymbol{R_{0}} = 0.84$.  Since $R_{0}<1$, by theorem 2.5 the disease free equilibrium point, $\boldsymbol{E_{0} = (98.99,0,0)}$  is globally asymptotically stable for the system $(2.18)-(2.20)$. The plots in  figure 1 for different initial conditions depict the global stability of the infection free equilibrium $E_{0}$ of the system $(2.18)-(2.20)$.\\

	\begin{figure}[hbt!]
		\begin{center}
			\includegraphics[width=5in, height=2.8in, angle=0]{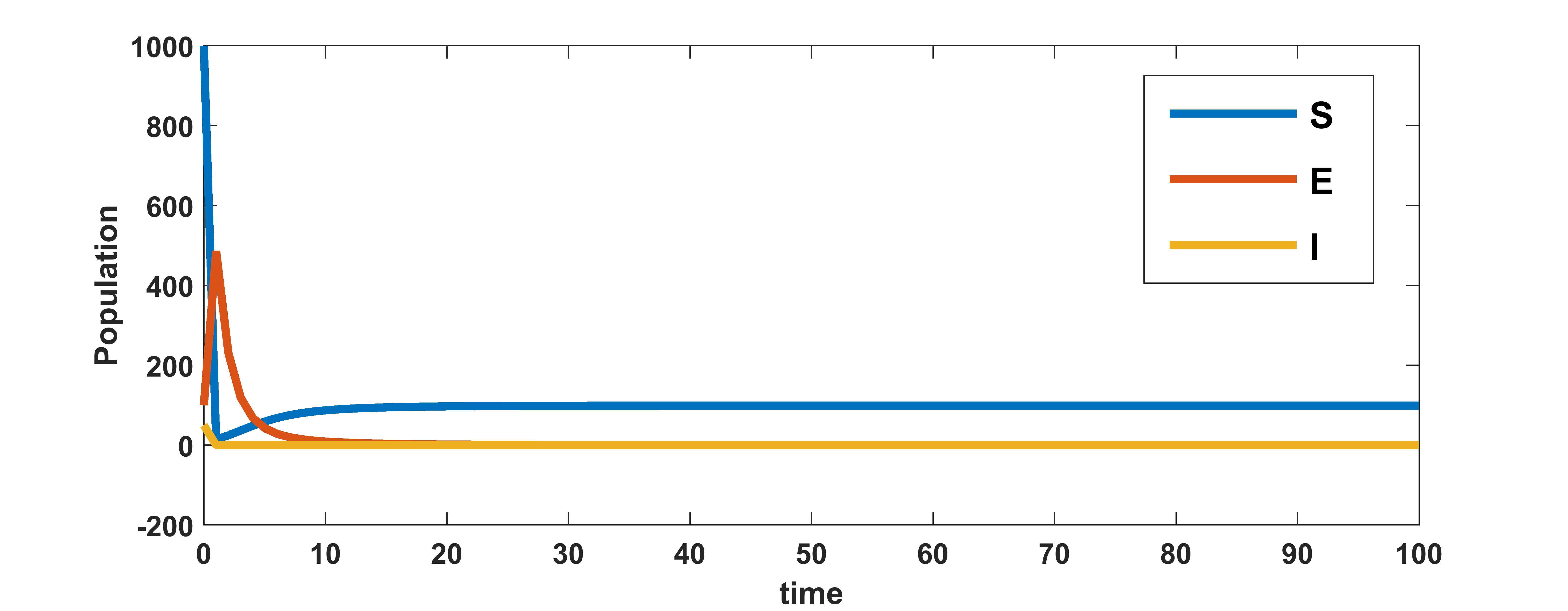}
			\includegraphics[width=5in, height=2.8in, angle=0]{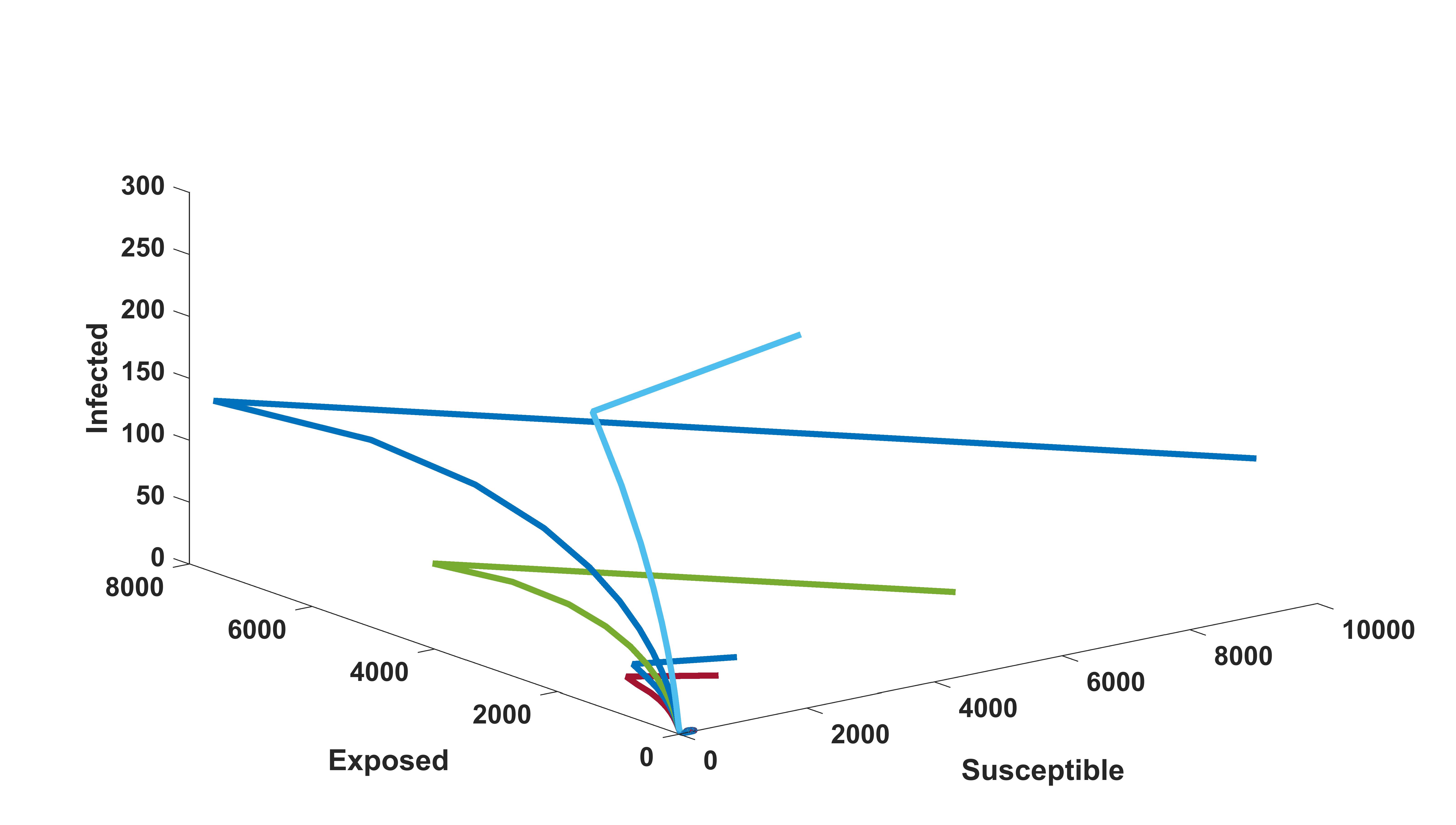}
			\caption{Figure depicting the  stability of $E_{0}$  of the system $(2.18)-(2.20)$. The first frame in the figure depicts the local asymptotic stability and the second frame depicts the global stability. }
			\label{sen_mu1_2}
		\end{center}
	\end{figure}

We know that the infected equilibrium $E_{1}$ exists only if $R_{0}>1$   and it  is also locally asymptotically stable whenever $R_{0}>1$. For the parameter values in the following table 3 the value of $R_{0}$ was calculated to be 135.7936 and  $E_{1}$ to be $(8.4687, 316.8, 165.03)$. Figure 2 demonstrates that $E_{1}$ is locally asymptotically stable whenever $R_{0} > 1$. In figure 3 the occurrence of the trans-critical bifurcation at $R_{0}=1$ is shown. 

\begin{figure}[hbt!]
\centering
	\includegraphics[height = 6cm, width = 15.5cm]{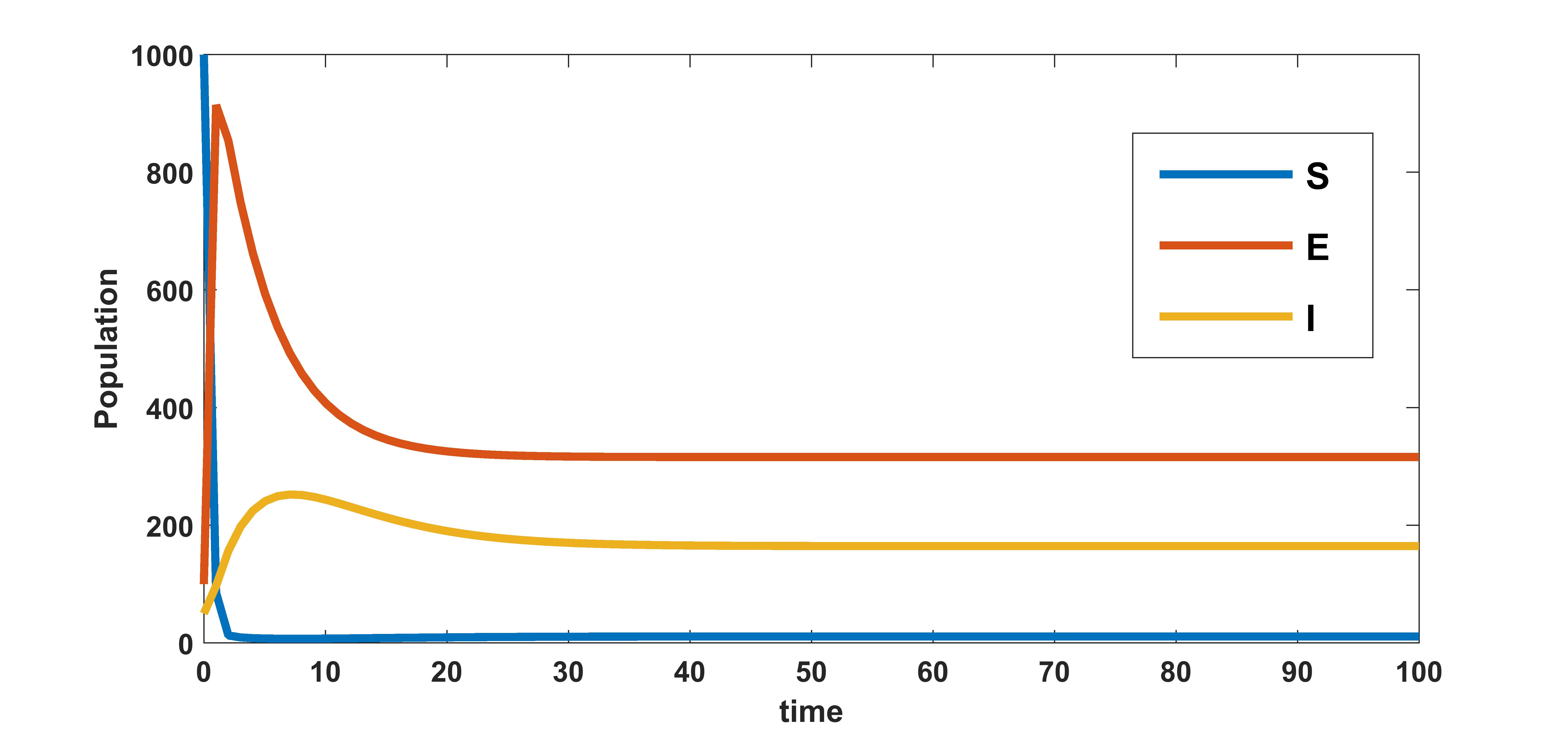}
	\caption{\label{first}Figure depicting the local asymptotic stability of $E_{1}$ of the system $(2.18)-(2.20)$ starting from the  initial state  $(S_{0}, E_{0}, I_{0}) = (1000, 100, 50).$ }
	\vspace{6mm}
\end{figure}

\begin{figure}[hbt!]
\centering
	\includegraphics[height = 6cm, width = 15.5cm]{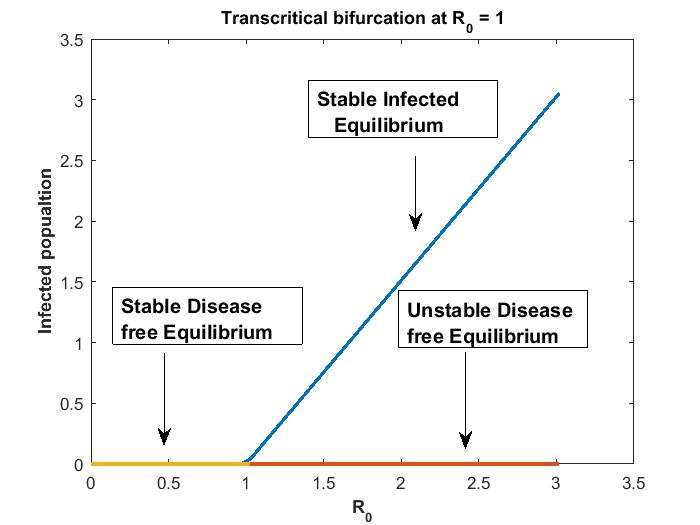}
	\caption{Figure depicting the trans-critical bifurcation  exhibited by the system $(2.18)-(2.20)$  at $R_0 = 1.$ The change in the stability of the equilibria with variation in $R_0$ can be  observed. }
\label{first}	
\end{figure}

\subsubsection{\textbf{Heat Plots}}
  Here we vary two parameters of the model $(2.18)-(2.20)$ at a time in a certain interval and plot the
the value of $R_0$ as heat plots. Heat plot has two different colours: one corresponding to the region with $R_0 <1$ and the other corresponding to $R_0>1$. This plot helps to find the combination of parameter values for which $R_0 <1$  and $R_0 > 1$. The blue colour region in these plots corresponds to the region where $R_0 < 1$ and
 therefore, from theorem 2.5, the disease free equilibrium is globally stable in this region. The other region with green colour corresponds to the region where $R_0 > 1$ and  In these region, the disease free equilibrium is unstable  and there exists a unique infected equilibrium point whose stability is determined using theorem 2.6. The stability of infected equilibria is determined using theorem 2.6. 
 
 \newpage
 \textbf{Parameter $\mu$ and $d$} \\
 Heat plot varying $\mu$ in the interval $(0.5, 0.94)$ and $d$ in the interval $(0.001, 0.44)$ is given in figure 4
 \begin{figure}[hbt!]
\centering
	\includegraphics[height = 8cm, width = 14cm]{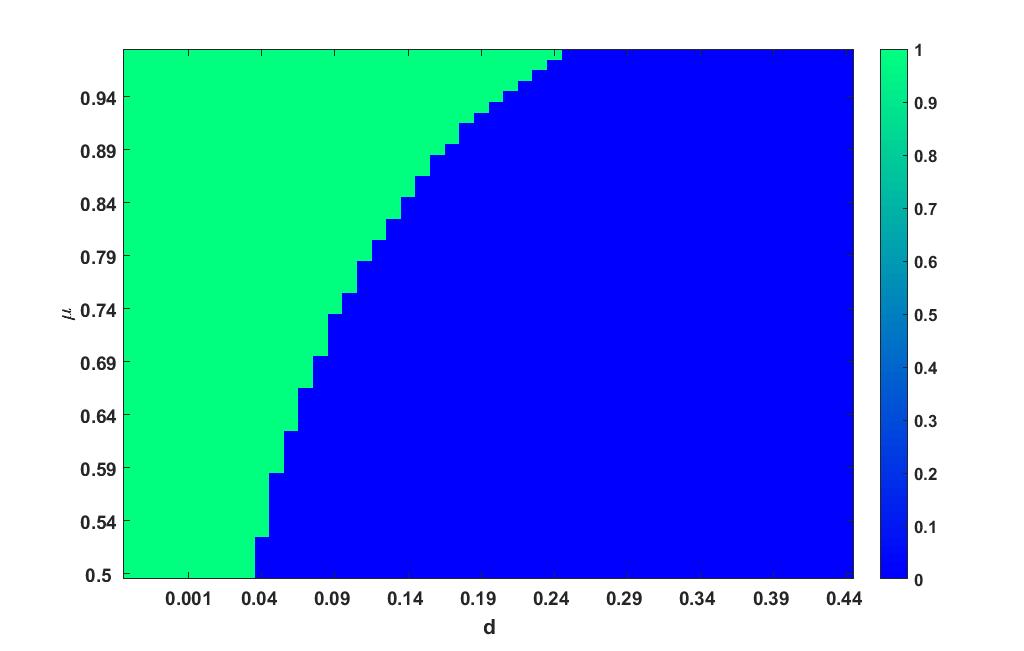}
	\caption{Heat plots for parameters $\mu$ and $d$}
\label{first}	
\end{figure}

  \textbf{Parameter $\pi$ and $d$} \\
 Heat plot varying $\pi$ in the interval $(0.0005, 0.005)$ and $d$ in the interval $(0.05, 0.5)$ is given in figure 5.

 \begin{figure}[hbt!]
\centering
	\includegraphics[height = 8cm, width = 14cm]{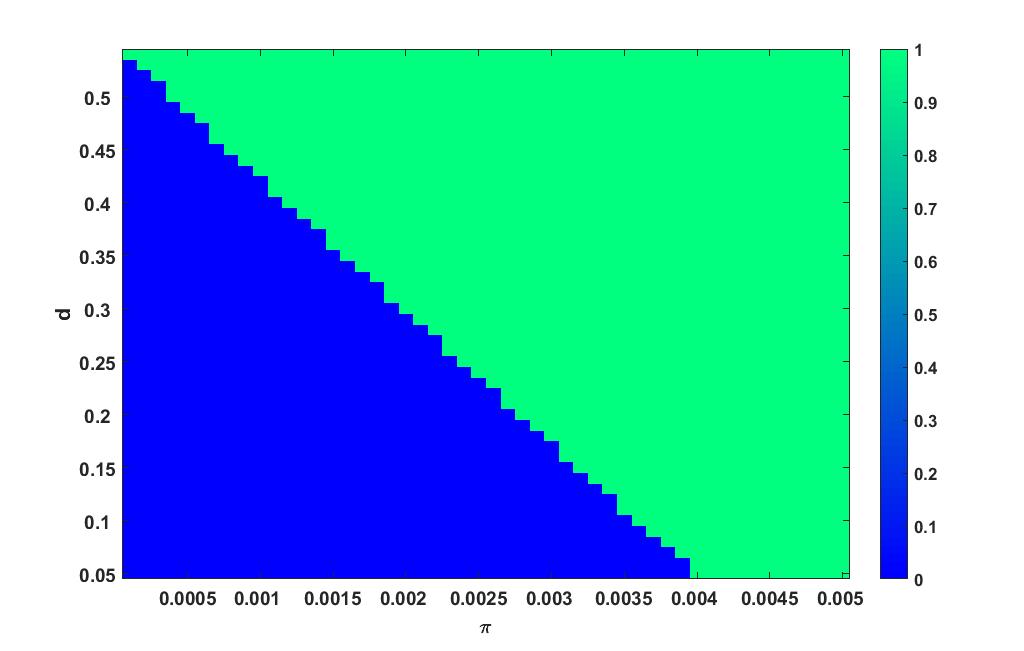}
	\caption{Heat plots for parameters $\pi$ and $d$}
\label{first}	
\end{figure}

   \textbf{Parameter $\mu$ and $\pi$} \\
 Heat plot varying $\pi$ in the interval $(0.004, 0.049)$ and $\mu$ in the interval $(1.05, 1.5)$ is given in figure 6.

 \begin{figure}[hbt!]
\centering
	\includegraphics[height = 8cm, width = 14cm]{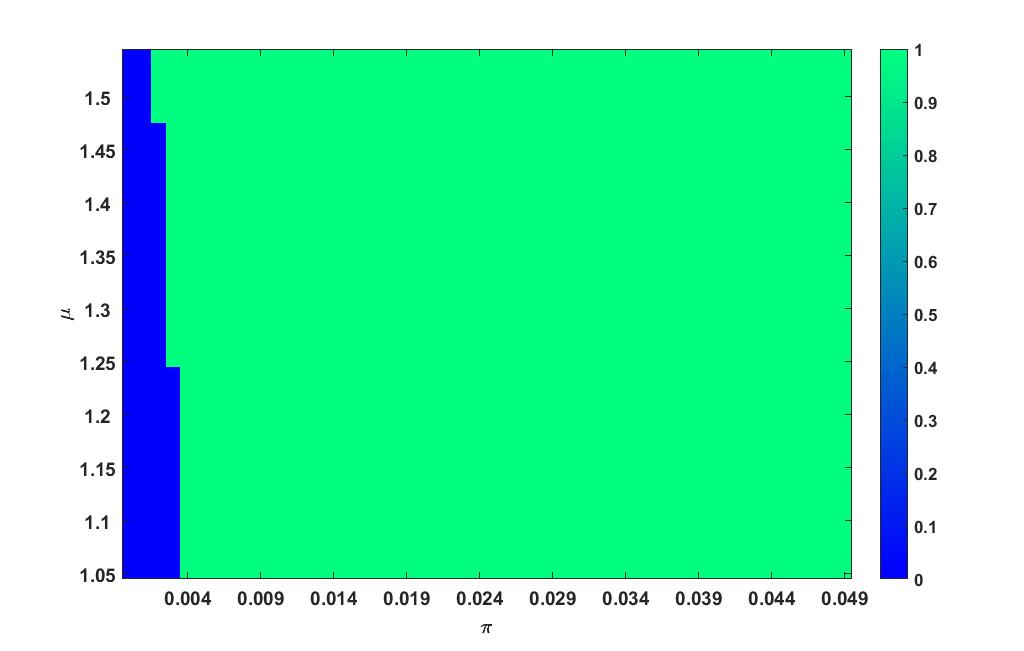}
	\caption{Heat plots for parameters $\mu$ and $\pi$}
\label{first}	
\end{figure}

  \textbf{Parameter $\beta$ and $d$} \\
 Heat plot varying $\beta$ in the interval $(0.0005, 0.005)$ and $d$ in the interval $(0.05, 0.5)$ is given in figure 7.
 \begin{figure}[hbt!]
\centering
	\includegraphics[height = 8cm, width = 14cm]{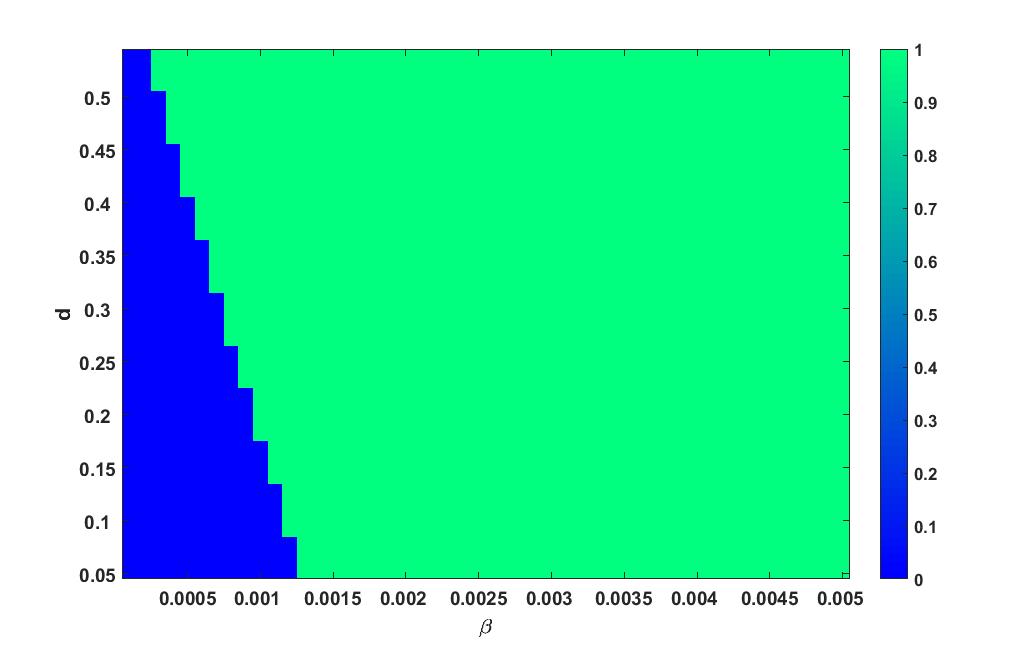}
	\caption{Heat plots for parameters $\beta$ and $d$}
\label{first}	
\end{figure}

 \subsubsection{\textbf{Influence of within-host sub-model parameters 
on the between-host sub-model}}

In this section we study the influence of the within-host sub-model parameters on the between-host sub-model variables. The reduced multi-scale model $(2.18)-(2.20)$ is categorized as a nested multi-scale model according to the categorization of multi-scale models infectious disease systems \cite{garira2017complete}. Therefore, the  multi-scale model $(2.18)-(2.20)$ is unidirectionally coupled. In this only the within-host scale sub-model influences the between-host scale sub-model without any reciprocal
feedback. Here we illustrate the influence of the key within-host scale sub-model parameters such as $\alpha, y$ and $x$ on the between-host scale sub-model variable $I$. At the within-host scale sub-model, $\alpha, y$ and $x$ describes the production of virus by infected cells, the clearance of the infected cells by the immune system, and the clearance
rate of free virus particles by the immune system. The parameter values for the simulation are taken from table 3. \\

In figure 8, the effect of variation
of burst rate of virus on infected population is plotted. The infected population $I(t)$ is plotted for three different values of burst rate $\alpha=0.24$, $\alpha=0.5$, and $\alpha=0.7$. These numerical illustration show that the between-host scale variable $I(t)$ is influenced by the within-host scale parameter $\alpha$. We see from figure 8 that as
 the infected cell burst size increases, transmission of the disease in the community also increases. Therefore, antiviral drugs such as remdesivir, arbidol, and HCQ that reduce the average infected cell viral production rate at within-host scale will possibly reduce the transmission of COVID-19 disease  at between-host scale.

	\begin{figure}[hbt!]
		\begin{center}
			\includegraphics[width=5in, height=2.8in, angle=0]{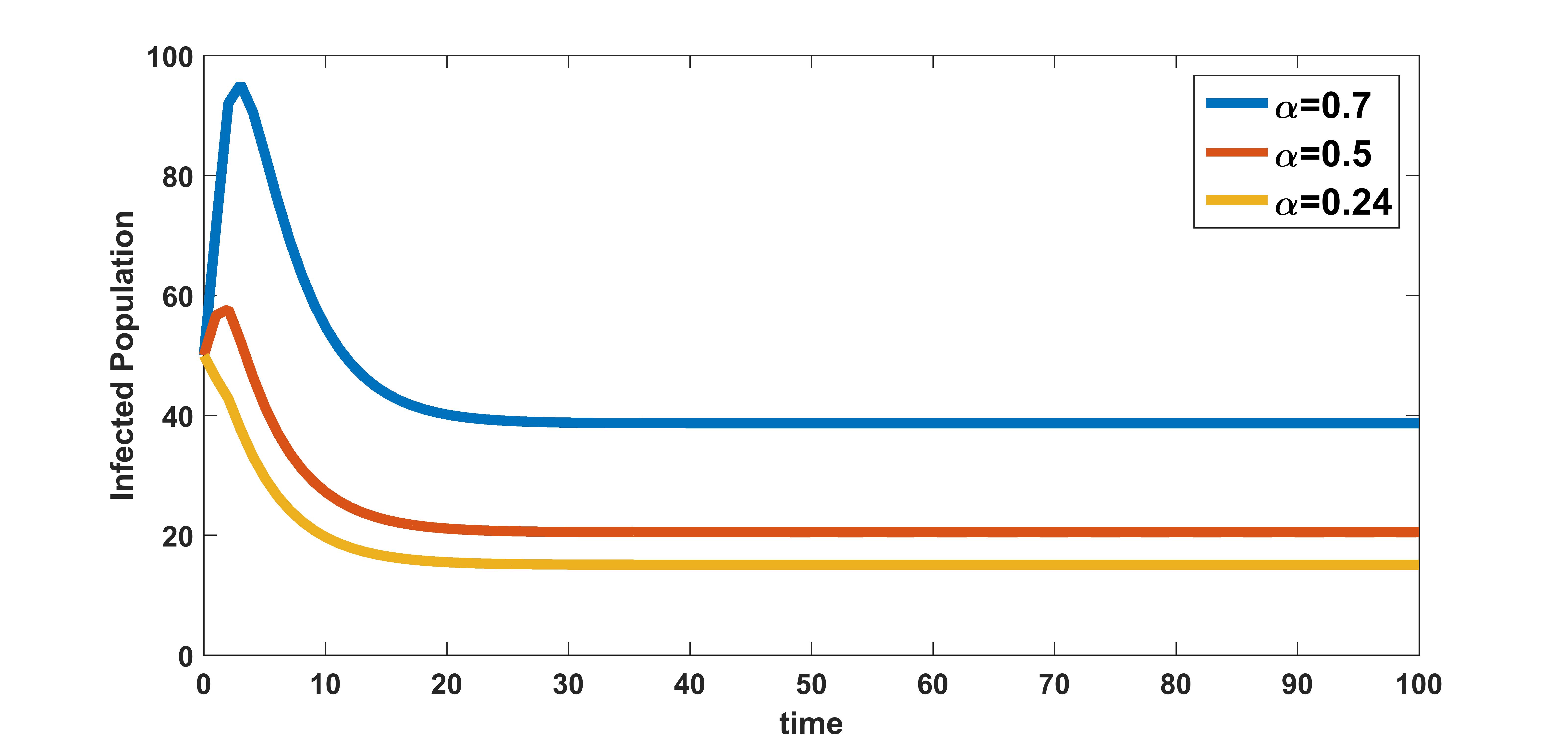}
			
			\caption{ Effect of variation
of burst rate of virus on infected population}
			\label{sen_mu1_2}
		\end{center}
	\end{figure}

In figure 9, the effect of variation
of rate of clearance of infected cells by the immune system on infected population is plotted. The infected population $I(t)$ is plotted for three different values  $x=0.5$, $x=0.795$, and $x=0.85$. These numerical illustration show that the between-host scale variable $I(t)$ is influenced by the within-host scale parameter $x$. We see from figure 9 that as
 rate of clearance of infected cells by the immune system increases, the infection at population level decreases. Therefore, drugs
that kill infected cells could have community level benefits of reducing COVID-19 transmission.

	\begin{figure}[hbt!]
		\begin{center}
		
			\includegraphics[width=5in, height=2.8in, angle=0]{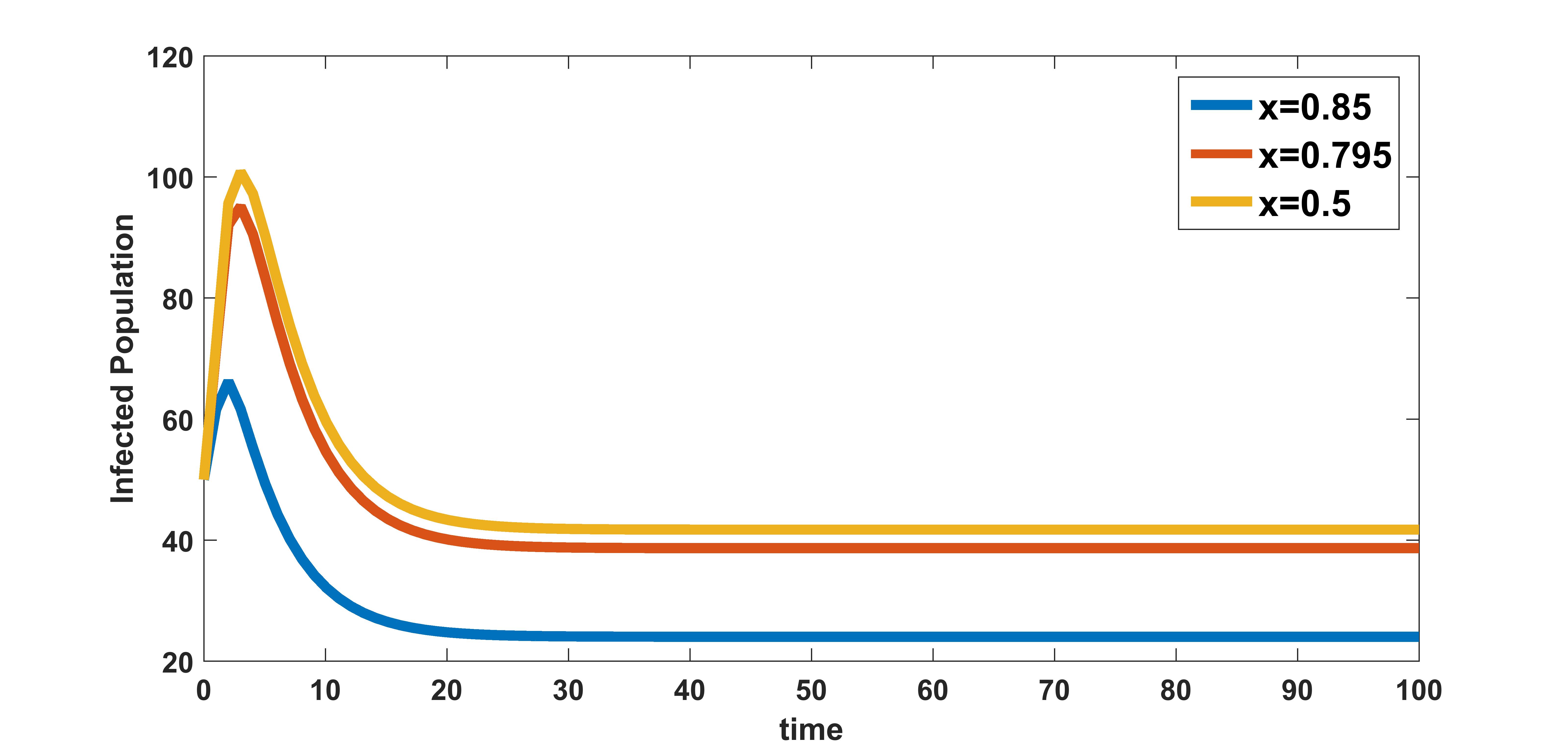}
			\caption{ Effect of variation
of rate of clearance of infected cells by the immune system on infected population}
			\label{sen_mu1_2}
		\end{center}
	\end{figure}

In figure 10, the effect of variation
of rate of clearance of virus particles by the immune system on infected population is plotted. The infected population $I(t)$ is plotted for three different values  $y=0.56$, $y=0.7$, and $y=0.8$. These numerical illustration show that the between-host scale variable $I(t)$ is influenced by the within-host scale parameter $y$. We see from figure 10 that as
 rate of clearance of virus by the immune system increases, the infection in the community decreases. Therefore, in addition to the benefits at individual level, treatments that increase the clearance rate of free virus particles in an infected individual could also have potential benefits in reducing the transmission at between-host level.

\begin{figure}[hbt!]
		\begin{center}
		
			\includegraphics[width=5in, height=2.8in, angle=0]{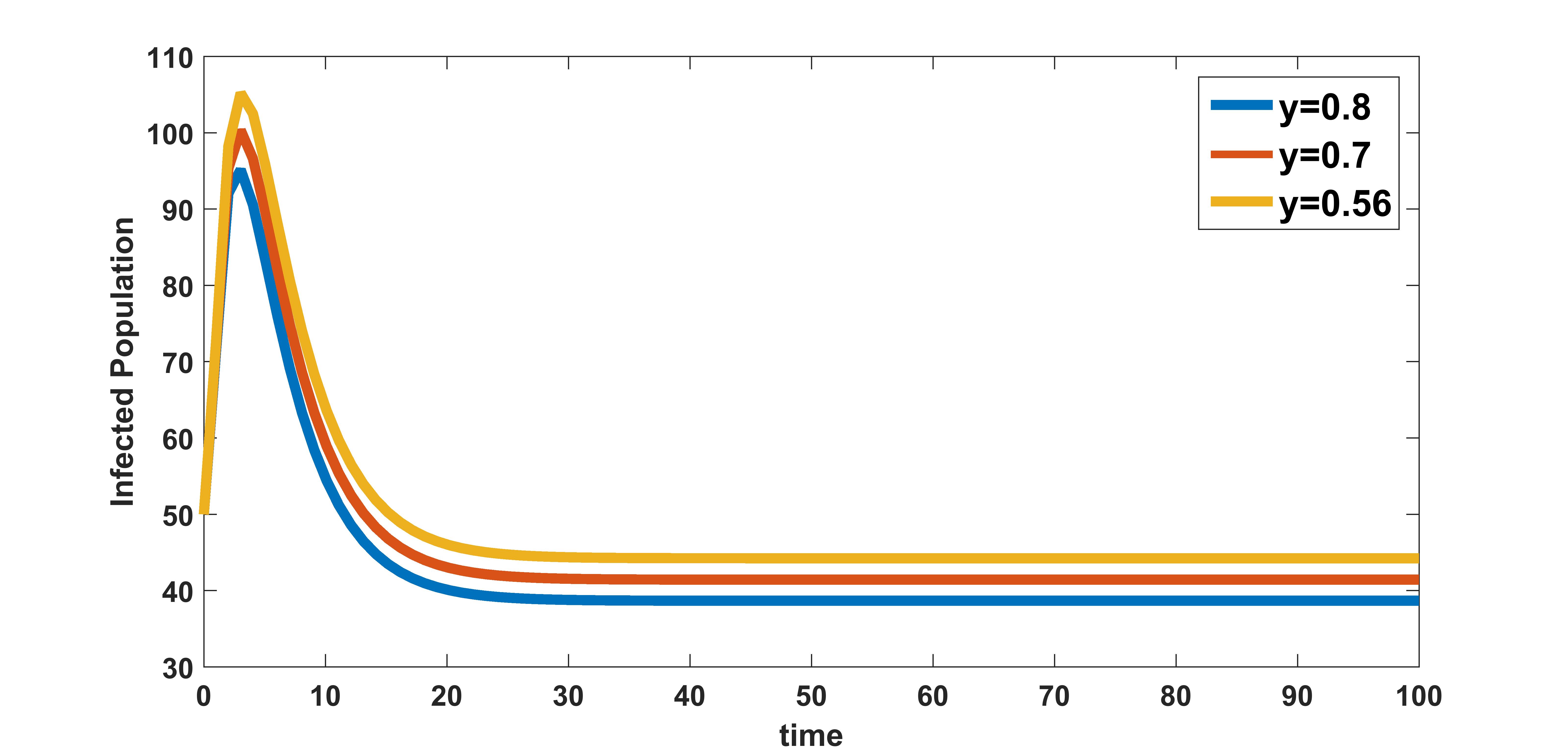}
			\caption{ Effect of variation
of rate of clearance of virus by the immune system on infected population}
			\label{sen_mu1_2}
		\end{center}
	\end{figure}

\section{Comparative Effectiveness of Health Interventions}	
		
In this section, in similar lines to the study done in \cite{garira2020development, prakash2020control} we extend the work on the multi-scale model discussed in Section 2 by including health care interventions. In the context of infectious diseases, disease dynamics can cause a large difference between the performance of a health intervention at the individual level (within-host scale), which is easy to determine, and its performance at the population level (between-host scale), which is difficult to determine \cite{garira2020development}. As a result, in situations where the effectiveness of a health intervention cannot be determined, health interventions with proven effectiveness may be recommended over those with potentially higher comparable effectiveness \cite{garira2020development}. We use three health care interventions that are as follows:\\

(1) Antiviral Drugs: 

(a) Drugs such as Remdesivir inhibit RNA-dependent RNA polymerase and drugs   Lopinavir and Ritonavir  inhibit the viral protease there by reducing  viral replication \cite{tu2020review}. Considering
these antiviral drugs are administered the burst rate of the within-host scale sub-model $\alpha$ now gets modified to $ \alpha(1-\epsilon)$ where $\epsilon$ is the efficacy of antiviral drugs and $0 < \epsilon < 1.$ 

(b) HCQ acts by preventing the SARS-CoV-2 virus from binding to cell membranes and also CQ/HCQ could  inhibit viral entry by acting as inhibitors of the biosynthesis of sialic acids, critical
actors of virus-cell ligand recognition \cite{roldan2020possible}. The administration of this  drug decreases the infection rate $k$ of susceptible
cells.This parameter gets modified to become $k(1-\gamma)$ where $\gamma$ is the efficacy of these drug and $0 < \gamma < 1$.

(2) Immunomodulators: Interferons are broad spectrum  antivirals, exhibiting both direct inhibitory effect on viral replication and supporting an immune response to clear virus infection \cite{wang2019global}. Due to the administration of immunomodulators such as INF immunity power of an individual increases as a result  of which the clearance rate of infected cells and virus particles by immune response increases. Therefore the parameter $x$ and $y$ now gets modified to $x(1+\delta)$ and $y(1+\delta)$ $\delta$ is the efficacy of immunomodulators and $0 < \gamma < 1$.

(3) Generalized social distancing: 
This intervention involves introducing measures such as restriction of mass gatherings, wearing of mask,
temporarily closing schools etc.  Assuming
that generalized social distancing is introduced to control COVID-19 epidemic, then the rate of contact with community $\beta$ gets modified to $\beta(1-\rho)$ where $\rho$ is the efficacy of generalized social distancing and $0 < \rho <1$.

Considering all the modifications in the parameters, the multi-scale model incorporating the effects of all the above health interventions  becomes

\begin{eqnarray}
 	\frac{dS}{dt}& =&  \Lambda \ - \beta(1-\rho) N_m S(t)I(t)  - \mu S(t)  \label{sec2equ1} \\
   	\frac{dE}{dt} &=& \beta N_m S(t)I(t) \ - (\mu+\pi+\gamma_1)E(t)    \label{sec2equ2}\\ 
   	\frac{dI}{dt} &=& \pi E(t)-(\mu+\gamma_2)I(t)-d N_m I(t)\label{sec2equ3} 
   \end{eqnarray}

where $N_m$ is modified $N_h$ given by
 \begin{equation}
     N_m = \frac{\alpha(1-\epsilon) \int_{d_1}^{d_2}U^* ds}{y(1+\delta)+\mu_v}
 \end{equation}

The effective reproduction number after incorporating the health interventions is given by,
$$R_{E}=\frac{\beta (1-\rho) N_m \Lambda \pi }{\mu(\mu+\pi+\gamma_1)(\mu +\gamma_2+d N_m)}$$

Here the comparative
effectiveness of the above three health interventions  are evaluated using $R_E$ as indicators of intervention effectiveness. In order to evaluate these we calculate the percentage reduction of $\mathcal{R}_{0}$  for single and multiple combination of the interventions at different efficacy levels such as (a) low efficacy of $0.3$, (b) medium efficacy of $0.6$, and (c)high efficacy of $0.9$. The different efficacy levels of the health interventions are chosen from \cite{prakash2020control}.

Percentage reduction of $\mathcal{R}_{0}$ is given by
$$ \bigg[ \frac{\mathcal{R}_{0} - \mathcal{R}_{E_j}}{\mathcal{R}_{0}} \bigg] \times 100$$
where $R_{Ej}$ means effective reproductive number when intervention/combination of interventions with efficacy/efficacies
$j$.

We now consider $8$ different combinations of these three health interventions corresponding to efficacy values 0.3, 0.6, and 0.9 obtained using the effective reproductive number $R_E$ as the indicator of intervention effectiveness. Then for each efficacy level, we rank the percentage reductions on $\mathcal{R}_{0}$ in ascending order from 1 to 8 corresponding to the different combinations of three health interventions considered in this study. The comparative effectiveness is calculated and measured on a scale from $1$ to $8$ with $1$ denoting the lowest comparative effectiveness and $8$ denoting the highest comparative effectiveness. In Table 5 the  abbreviations a) CEL stands for "Comparative Effectiveness at Low efficacy," which is 0.3, b) CEM stands for "Comparative Effectiveness at Medium efficacy," which is 0.6, c) CEH stands for "Comparative Effectiveness at High efficacy" which is 0.9.

\vspace{.25cm}

\begin{table}[htp!]
\caption{Comparative effectiveness for $\mathcal{R}_{0}$}
\begin{center}
\begin{tabular}{ | c | c | c | c | c | c | c | c |}
\hline
 \textbf{No.} & \textbf{Indicator} & \textbf{\%age} & \textbf{CEL} & \textbf{\%age} & \textbf{CEM} & \textbf{\%age} & \textbf{CEH} \\ 
 \hline \hline
 1 & $\mathcal{R}_{0}$ & 0  &  1  &  0  &  1  &  0  &  1 \\
 \hline \hline
 2 & $\mathcal{R_{E_{\rho}}}$ & 30  &  5 &  60  &  5  &  90  &  5 \\
 \hline
 3 & $\mathcal{R_{E_{\delta}}}$ & 0.106  &  3  &  0.24  &  2  &  0.45  &  4  \\ 
 \hline \hline
 4 & $\mathcal{R_{E_{\epsilon}}}$ & 0.075  &  2  &  0.26  &  3  &  0.38  &  3  \\ 
 \hline
 5 & $\mathcal{R_{E_{\rho \delta}}}$ & 30.07  &  7  &  60.12  &  7  &  90.45  &  7  \\
 \hline \hline
 6 & $\mathcal{R_{E_{\rho \epsilon}}}$ & 30.05  &  6  &  60.1  &  6  &  90.23  &  6  \\
 \hline \hline
 7 & $\mathcal{R_{E_{\epsilon \delta}}}$ & 0.22  &  4  &  0.35  &  4  &  0.67  &  2  \\
 \hline \hline
 8 & $\mathcal{R_{E_{\rho\delta \epsilon}}}$ & 30.16  &  8  &  60.34  &  8  &  90.5  &  8  \\ 
 
 \hline\hline
\end{tabular}

\label{cet1}
\end{center}
\end{table}

From table 5, we deduce the following results regarding comparative effectiveness of three different interventions considered: 

(a) When a single intervention strategy is implemented, intervention involving introducing measures such as restriction of mass gatherings, wearing of mask,
temporarily closing schools etc. show significant decrease in $R_0$ compared to the implementation of antiviral drugs and immunomodulators at all efficacy levels.   

(b) When considering two interventions, we observe that the generalized social distancing along with  Immunomodulators that boost the immune response would be highly effective in limiting the spread of infection in the community.

(c) A combined strategy involving treatment with anitviral drugs, immunomodulators, and generalized social distancing seems to perform the best among all the combinations considered at all efficacy levels.

From this section, we conclude that a combined strategy is the best strategy to contain the spread of infection. However, the implementation of a combined strategy may not always be cost-effective and feasible. In a single intervention strategy, generalized social distancing is shown to lower the value of $R_0$ better than individual use of antiviral drugs and immunomodulators. Therefore, in a situation where resources are limited and costly, interventions such as restricting mass gatherings, wearing masks, temporarily closing schools, etc. would be highly effective in containing infection and also incur less cost.

\section{Discussions and Conclusions}
\noindent
 
 COVID-19 is a contagious respiratory and vascular disease that has resulted in more than 195 million cases and 4.18 million deaths worldwide. It is  caused by severe acute respiratory syndrome coronavirus 2 (SARS-CoV-2). On 30 january it was declared as a public health emergency of international concern \cite{concern}. Mathematical models have proved to provide useful information about the dynamics of  the infectious diseases.  To understand the dynamics  of the COVID-19 disease several compartment models has been developed \cite{samui2020mathematical, ndairou2020mathematical, zeb2020mathematical,leontitsis2021seahir, wang2020four, dashtbali2021compartmental, zhao2020five, chen2020mathematical, biswas2020covid, sarkar2020modeling}. 

An in-depth understanding of the transmission of infectious disease systems requires
knowledge of the processes at the various scales of infectious diseases and how these scales
interact. In this study, we develop a nested multi-scale model for COVID-19 disease that integrates the within-host scale and the between-host scale sub-models. The transmission rate and COVID-19 induced death rate at between-host scale are assumed to be a linear function of viral load. Because of the difficulties of working at two different time scales, we approximate individual-level host infectiousness   by some surrogate measurable quantity called area under viral load and reduce the multi-scale model at two different times scales to a model with area under the viral load curve acting as a proxy of individual level host infectiousness. 

Initially, the well posedness of the reduced multi-scale model is discussed followed by the stability analysis of the equilibrium points admitted by the reduced multi-scale model. The disease free equilibrium point of the model is found to remain globally asymptotically stable whenever the value of basic reproduction number is less than unity. As the value of basic reproduction  number crossed unity a unique infected equilibrium point exists and remains stable depending on the sign of $(AB-C)$. The model is shown to undergo a forward (trans-critical) bifurcation at $R_0=1$. To predict the sensitivity of the model parameters  on $R_0$, elastic index that measures the
relative change of $R_0$ with respect to parameters is calculated for each parameters in the definition of $R_0$. The parameters $\beta$, $\pi$ and $\Lambda$ were found to be most sensitive towards $R_0$. The theoretical results are supported with numerical illustrations. To separate out the region of stability and instability of the equilibrium points, two parameter heat plot is done by varying the parameter values in certain range. The influence of the key within-host scale sub-model parameters such as $\alpha,y$ and $x$ on between-host scale are also numerically illustrated. It is found that the spread of infection in a community is influenced by the production of virus particles by an infected cells ($\alpha$), clearance rate of infected cell $(x)$ and clearance rate of virus particles by immune system $(y)$.

We also use the reduced multi-scale model developed to study the comparative effectiveness of the three health interventions (antiviral drugs, immunomodulators and generalized social distancing) for COVID-19 viral infection using $R_E$ as the indicator of intervention effectiveness. The result suggested that a combined strategy involving treatment with anitviral drugs, immunomodulators, and generalized social distancing would be the best strategy to limit the spread of infection in the community. However, implementing a combined strategy may not be always cost-effective or feasible. In a single intervention strategy, general social distancing has been shown to lower the value of $R_0$ better than individual use of antiviral drugs and immunomodulators. Therefore, in a situation where resources are limited and costly, interventions such as restricting mass gatherings, wearing masks, temporarily closing schools, etc., would be highly effective in containing infection and also more cost-effective. We believe that the
results presented in this study will help physicians, doctors and researchers in making informed decisions about COVID-19 disease prevention and treatment interventions.

{\flushleft{  \textbf{DEDICATION} }}\vspace{.25cm}

The authors from SSSIHL  dedicate this paper to the founder chancellor of SSSIHL, Bhagawan Sri Sathya Sai Baba. The first author dedicates this paper to his loving father Purna Chhetri and the second author to his loving elder brother D. A. C. Prakash who still lives in his heart. 
   \bibliographystyle{amsplain}
\bibliography{reference}

\end{document}